\documentclass[11pt, reqno]{amsart}

\usepackage[usenames,dvipsnames]{pstricks}
\usepackage{epsfig}

\usepackage{color}
\date{\today}

\usepackage[latin1]{inputenc}
\usepackage[T1]{fontenc}
\usepackage{amsmath, amsthm}
\usepackage[english]{babel}
\usepackage{amssymb}
\usepackage{latexsym}
\usepackage{fancyhdr}
\usepackage{graphicx}
\usepackage[all]{xy}

\theoremstyle{definition}
\newtheorem{definicion}{{\bf Definition}}[section]
\newtheorem{parrafo}[definicion]{}
\newtheorem{nota}[definicion]{{\bf Remark}}
\newtheorem{exemplo}[definicion]{{\bf Example}}

\newtheorem{teorema}[definicion]{{\bf Theorem}}
\newtheorem{prop}[definicion]{{\bf Proposition}}

\newtheorem{corol}[definicion]{{\bf Corollary}}

\newenvironment{dem}%
{{\noindent{\bf Proof: }}\newline }%
{$\hfill\Box$\vspace{0.25cm}} %

\begin{document}

\begin{center}
 {\huge{\bf  Weak Crossed Biproducts and Weak }}

\vspace{0.25cm}

 {\huge{\bf Projections}}
\end{center}

\ \\

{\bf  J.M. Fern\'andez Vilaboa$^{1}$, R. Gonz\'{a}lez
Rodr\'{\i}guez$^{2,\ast}$ and A.B. Rodr\'{\i}guez Raposo$^{3}$}

\ \\
\hspace{-0,5cm}$1$ Departamento de \'Alxebra, Universidad de
Santiago de Compostela,  E-15771 Santiago de Compostela, Spain
(e-mail: josemanuel.fernandez@usc.es)
\ \\
\hspace{-0,5cm}$2$ Departamento de Matem\'{a}tica Aplicada II,
Universidad de Vigo, Campus Universitario Lagoas-Marcosende,
E-36310 Vigo, Spain (e-mail: rgon@dma.uvigo.es)
\\
\hspace{-0,5cm}$3$ Departamento de Computaci\'on, Universidade da Coru\~{n}a,  E-15071 A Coru\~{n}a, Spain
(e-mail: arodriguezraposo@udc.es)
\ \\
\hspace{-0,5cm}$\ast$ Corresponding author.

\begin{center}
{\bf Abstract}
\end{center}
{\small In this paper we present the general theory and universal properties of weak crossed biproducts. We prove that every weak projection of weak bialgebras induces one of these weak crossed structures. Finally, we compute explicitly the weak crossed biproduct associated to a groupoid that admits an exact factorization.}

\vspace{0.5cm}



\section*{Introduction}

With the recent arise of some kinds of weak algebraic structures, and more specifically of weak Hopf algebras, several important concepts in the theory of Hopf algebras had to be adapted to this new context. In particular, none of the crossed product type structures in the Hopf algebra setting  is suitable for the weak setting, because of the lack of some properties of the units and the counits. The purpose of this paper is to obtain a general and appropriate framework  to develop all types of crossed products, crossed coproducts and crossed biproducts (i.e., weak bialgebras given as crossed products and coproducts) in the weak setting, that moreover generalizes the classical crossed product algebraic structures like crossed product algebras, bicrossproducts or cross product bialgebras among others. In particular, we are interested in the crossed biproduct structure that arises from a weak bialgebra with a weak projection.

Crossed products appear as a generalization of semi-direct products of groups to the context of Hopf algebras \cite{Molnar,bcm}, and were also used to study cleft extensions and Galois extensions of Hopf algebras \cite{blat-susan, doi3}. In \cite{tb-crpr} Brzezi\'nski gave a categorical approach that generalizes several types of crossed products, even the ones given for braided Hopf algebras by Majid \cite{maj2}.  As Hopf algebras have a coalgebra structure, it is easy to obtain a theory of crossed coproducts in similar terms. Crossed products and crossed coproducts can be combined to obtain a Hopf algebra structure given in such terms. The first examples of crossed biproducts  arise in the study of Hopf algebras related to matched pairs of groups \cite{TAK} and the study of cohomology and exact sequences \cite{SING}. Those first examples suffered further generalizations to more and more general contexts, that are summarized using a categorical approach given by Bespalov and Drabant \cite{bes-drab1, bes-drab2}, that moreover, when it is only considered in its algebra version, recovers classical crossed products  \cite{bcm,doi3, maj1, maj99, maj2}.  Bespalov and Drabant define a general bialgebra structure, a cross product bialgebra, by means of a crossed product without the explicict use of a cocycle and a crossed coproduct without the explicit use of a cycle. They furthermore obtain some universal properties that characterize cross product bialgebras, and hence its categorical meaning is completely understood. Their general framework unifies all known non weak crossed biproduct type constructions in a single setting, like the Drinfeld's quantum double,  Radford's 4-parameter Hopf algebra, the quantum Poincar\'e group,  the quantum Weyl group, Lusztig's construction of the quantum enveloping algebra, the affine quantum groups $U_{q}($\^{g}$)$, the Connes-Moscovici Hopf algebra, etc. Moreover, it provides more new and interesting examples of cross product bialgebras in braided contexts.

A relevant example of a crossed biproduct arises from Hopf algebras with a projection. Consider a bialgebra $D$ and a Hopf algebra $B$ such that $f:B\rightarrow D$ and $g:D\rightarrow B$ are morphisms of bialgebras that satisfy $g\circ f = id_B$. In \cite{RAD} Radford shows that $D$ can be recovered as a crossed biproduct whose algebra and coalgebra structures  are given by the smash product and the cosmash coproduct respectively (see \cite{maj2} for the braided setting). If the hypothesis on $g$ are relaxed in an appropriate way we still obtain the algebra $D$ as a crossed biproduct structure, and although the algebra is not an smash product, the coalgebra is still an smash coproduct. Suppose that $g$ is no longer a morphism of bialgebras but it is a morphism of coalgebras that is also a morphism of right $B$-modules, that is, $g\circ\mu_D \circ (D\otimes f) = \mu_B\circ (g\otimes B)$ for $\mu_B$ and $\mu_D$ the product on $B$ and $D$ respectively. In this case we say that $D$ has a weak projection onto $B$. For example, quantized enveloping algebras have a weak projection, and also  finite dimensional pointed Hopf algebras have a weak projection onto its coradical. Working in a category of vector spaces, Schauenburg  analizes this situation in \cite{SCH2} and obtains that, if $D$ and $B$ are bialgebras, $D$ is obtained as a bicrossproduct type structure. As expected, the coalgebra structure is a cosmash coproduct. The braided version was studied by Ardizzoni, Menini and Stefan in \cite{ardi}. In this case they use the categorical approach  given by Bespalov and Drabant in \cite{bes-drab1, bes-drab2}.

Although the notions and the results cited below are useful tools to study relevant examples in the Hopf algebra world,  they are not general enough to solve analogous problems in the weak Hopf algebra case. Weak Hopf algebras (or quantum groupoids) were recently introduced by B\"ohm, Nill and Szlach\'anyi in \cite{bohm1}. Recall that although a weak Hopf algebra has a multiplication and a comultiplication that are compatible like in the non weak case, the unit is not comultiplicative (or, equivalently, the counit is not multiplicative). However, there still exists a relation between the unit and the comultiplication of the weak Hopf algebra (and of the counit and the multiplication), that reveal the existence of some idempotent morphisms which in the non weak context become trivial. As a consequence, when one tries to extend crossed product constructions to the weak setting one finds that idempotent morphisms arise in a very natural way.  The use of idempotents combined with the ideas  in \cite{tb-crpr} are the key to introduce a general theory of weak crossed products.  In \cite{nmra4} we define a
product on $A\otimes V$, for an algebra $A$
and  an object $V$  both living in a strict monoidal category $\mathcal C$ where every idempotent splits. To obtain that product we consider two morphisms
$\psi_{V}^{A}:V\otimes A\rightarrow A\otimes V$ and
$\sigma_{V}^{A}:V\otimes V\rightarrow A\otimes V$ that must satisfy some  twisted-like and cocycle-like conditions. Associated to these morphisms we define an idempotent $\nabla_{A\otimes V}:A\otimes V\rightarrow A\otimes V$, that becomes the identity in the non weak case. The image of this idempotent inherits the associative product from $A\otimes V$. In order to define a unit for $Im (\nabla_{A\otimes V})$, and hence to obtain an algebra structure, we require the existence of a morphism $\eta_V:K\rightarrow V$ satisfying certain properties. Unfortunately, this condition turned out to be restrictive, even in some natural examples as the crossed product induced by a weak cleft extension (see \cite{nmra1}). Hence in \cite{mra-preunit} we changed this  hypothesis by the existence of a preunit $\nu:K\rightarrow A\otimes V$. A preunit on a non unitary associative algebra $B$ is a morphism $\nu:K\rightarrow B$ that, very roughly speaking, is commutative and idempotent with respect to the product in $B$. Due to these particularities it is possible to define a natural idempotent $\nabla:B\rightarrow B$, that, in our case, is the morphism $\nabla_{A\otimes V}$. Moreover, the preunit induces a unit for $Im(\nabla_{A\otimes V})$, so it is an algebra as in the case studied in \cite{nmra4}. Using, like in \cite{nmra4}, morphisms $\psi_V^A$ and $\sigma_V^A$ it is possible to characterize weak crossed products with a preunit as products on $A\otimes V$ that are morphisms of left $A$-modules and that have a preunit. Finally, it is convenient to observe that, if the preunit is a unit, the idempotent becomes the identity and we recover the classical examples of the non weak setting. The theory presented in \cite{nmra4, mra-preunit} contains as a particular case the one developed by Brzezi\'nski in \cite{tb-crpr}. There are many other examples of this theory like the weak smash product given by Caenepeel and De Groot in \cite{caengroot}, the theory of wreath products that we can find in \cite{LS} and the weak crossed products
for weak bialgebras given in \cite{ana1}. Recently, G. B\"ohm showed in \cite{bohm} that a monad in the weak version of the Lack and Street's 2-category of monads in a 2-category is identical to a crossed product system in the sense of \cite{nmra4}.

Now the next natural step is to explore bicrossed product type structures in the weak context, that is, to combine the theory of crossed products given in \cite{nmra4, mra-preunit} with its dual and obtain a general framework to deal with weak bialgebras coming from crossed products and crossed coproducts. A first example of the existence of a weak crossed biproduct appears in \cite{NikaRamon4}, where Radford's theory for Hopf algebras with projection is given for weak Hopf algebras. Again, as expected, the crossed product structure and the crossed coproduct structure obtained are smash-like ones, but in a weak context. Based in this example and in the work by Schauenburg \cite{SCH2} and Ardizzoni, Menini and Stefan \cite{ardi}, it seems natural to obtain also a weak  crossed biproduct structure that arises from a weak projection of weak bialgebras. We find an example of this situation in the theory of (finite) groupoids. Using the terminology and the results introduced by Mackenzie in \cite{MacK}, an exact factorization of a groupoid $G$ is a kind of semidirect product of two subgroupoids $V$ and $H$ that must satisfy certain conditions. Equivalently, $G$ is a vacant double groupoid or $V$ and $H$ are a matched pair of groupoids. These kind of groupoids were recently used by Andruskiewitsch and Natale \cite{AN, AN1} to obtain a new class of weak Hopf algebras. If we consider the weak Hopf algebras associated to the groupoid algebras $RG$, $RV$ and $RH$, for $R$ a commutative ring with unit, it turns out that $RG$ has a strictly weak projection onto $RV$, that is, the morphism between $RG$ and $RV$ is of right $RV$ modules but  it is not of algebras.

Inspired by the examples mentioned above, this paper is devoted to obtain a general theory of weak crossed biproducts. In the first section we present some basic facts about weak crossed products and coproducts, using the previous theory introduced in \cite{nmra4,mra-preunit}. Moreover we find that weak crossed products and coproducts are universal constructions that generalize results of Bespalov and Drabant \cite{bes-drab2}. In section 2 we use the results developed in section 1 to obtain a weak crossed biproduct structure on $A\otimes C$, for an algebra $A$ and a coalgebra $C$ living in a braided monoidal category. We define a weak crossed biproduct as a weak bialgebra whose algebra structure comes from a weak crossed product on $A\otimes C$ and whose coalgebra structure comes from a weak crossed coproduct on $A\otimes C$. These two structures are glued together by considering that the respective idempotent morphisms are the same. In the main result of this section, Theorem \ref{Teo-biproduct}, we find, using an universal construction, which conditions must satisfy a weak bialgebra to be obtained as a weak crossed biproduct. Section 3 is devoted to the study of weak bialgebras with a weak projection in a braided monoidal category. Here we obtain that a weak bialgebra with a weak projection satisfies the universal properties given in Theorem \ref{Teo-biproduct} and, as a consequence, it is recovered as a weak crossed biproduct whose coproduct is a cosmash like construction. The mentioned result recovers, as a particular instances, the ones given by Schauenburg \cite{SCH2} and by Ardizzoni, Menini and Stefan \cite{ardi}. Finally, in this last section we compute explicitly the weak crossed biproduct associated to a groupoid that admits an exact factorization.

Throughout the paper  ${\mathcal C}$ denotes a strict
 monoidal category with tensor product $\otimes$ and base object $K$. Given objects
$A$, $B$, $D$ and a morphism
$f:B\rightarrow D$, we write $A\otimes f$ for $id_{A}\otimes f$
and $f\otimes A$ for $f\otimes id_{A}$. Also we assume  that all
idempotent splits, i.e., for every morphism
$\nabla_{Y}:Y\rightarrow Y$, such that $\nabla_{Y}=\nabla_{Y}\circ
\nabla_{Y}$, there exist an object $Z$ and morphisms
$i_{Y}:Z\rightarrow Y$ and $p_{Y}:Y\rightarrow Z$ satisfying
$\nabla_{Y}=i_{Y}\circ p_{Y}$ and $p_{Y}\circ i_{Y}=id_{Z}$.

As for prerequisites, the reader is expected to  be familiar with
the notions of algebra (monoid), coalgebra (comonoid), module and comodule in the
monoidal setting.  Given an algebra  $A$ and a coalgebra
$C$, we let $\eta_{A}:K\rightarrow A$, $\mu_{A}:A\otimes A
\rightarrow A$, $\varepsilon_{D}: D\rightarrow K$, and $\delta_{D}
:D\rightarrow D\otimes D$ denote the unity, the product, the
counity, and the coproduct respectively. Also, if $A$, $B$ are algebras,
$f:A\rightarrow B$ is an algebra morphism if $f\circ \eta_{A}=\eta_{B}$ and $f\circ
\mu_{A}=\mu_{B}\circ (f\otimes f)$. In a dual form we have the notion of coalgebra
morphism and, if $A$ is an algebra and $D$ is a coalgebra, for two morphism
$f,g:D\rightarrow A$, the symbol $\wedge$ denotes the usual convolution product in the category ${\mathcal C}$, i.e., $f\wedge g=\mu_{A}\circ (f\otimes g)\circ \delta_{D}.$

In this paper, for an algebra $A$ and a coalgebra $C$, the triple $(A,C, \psi_{RR})$ denotes a right-right weak entwining structure, and therefore the following identies hold:
\begin{equation}
\label{a)}
\psi_{RR} \circ (C\otimes \mu_A ) = (\mu_A \otimes C)\circ
(A\otimes \psi_{RR} )\circ (\psi_{RR} \otimes A),
\end{equation}
\begin{equation}
\label{b)} (A\otimes \delta_C)\circ \psi_{RR} = (\psi_{RR} \otimes C)\circ
(C\otimes \psi_{RR} )\circ (\delta_C\otimes A),
\end{equation}
\begin{equation}
\label{c)}
\psi_{RR} \circ (C\otimes \eta_A) = (e_{RR}\otimes C)\circ
\delta_C,
\end{equation}
\begin{equation}
\label{d)}
(A\otimes \varepsilon_C)\circ \psi_{RR} = \mu_A\circ
(e_{RR}\otimes A),
\end{equation}
where $e_{RR} = (A\otimes \varepsilon_C)\circ \psi_{RR} \circ (C\otimes
\eta_A)$. We denote by ${\mathcal
M}_A^C(\psi_{RR})$ the category of weak entwined modules, i.e., an object $M$
in ${\mathcal C}$ together with two morphisms $\phi_M:M\otimes
A\rightarrow A$ and $\rho_M:M\rightarrow M\otimes C$ such that
$(M, \phi_M)$ is a right $A$-module, $(M, \rho_M)$ is a right
$C$-comodule and
\begin{equation}
\label{m)}
\rho_M\circ \phi_M = (\phi_M\otimes C)\circ (M\otimes \psi_{RR} )
\circ (\rho_M\otimes A).
\end{equation}
The morphisms in ${\mathcal
M}_A^C(\psi_{RR})$ are the obvious, i.e., morphisms of right $A$-modules and right $C$-comodules.

\section{Weak crossed products and coproducts}

In this section we develop the general theory of weak crossed products in a monoidal cateory $\mathcal C$. In order to make this paper as self-contained as possible, we recall some results given in \cite{mra-preunit} that will be used to obtain further characterizations of weak crossed products and weak crossed coproducts.

Let $A$ be an algebra and $V$ be an object in
${\mathcal C}$. Suppose that there exists a morphism
$\psi_{V}^{A}:V\otimes A\rightarrow A\otimes V$  such that the following
equality holds
\begin{equation}\label{wmeas-wcp}
(\mu_A\otimes V)\circ (A\otimes \psi_{V}^{A})\circ
(\psi_{V}^{A}\otimes A) = \psi_{V}^{A}\circ (V\otimes \mu_A).
\end{equation}
 As a consequence of (\ref{wmeas-wcp}), the morphism $\nabla_{A\otimes V}:A\otimes V\rightarrow
A\otimes V$ defined by
\begin{equation}\label{idem-wcp}
\nabla_{A\otimes V} = (\mu_A\otimes V)\circ(A\otimes
\psi_{V}^{A})\circ (A\otimes V\otimes \eta_A)
\end{equation}
is  idempotent. Moreover, $\nabla_{A\otimes V}$
satisfies that $$\nabla_{A\otimes V}\circ (\mu_A\otimes V) =
(\mu_A\otimes V)\circ (A\otimes \nabla_{A\otimes V}),$$ that is,
$\nabla_{A\otimes V}$ is a left $A$-module morphism (see Lemma 3.1 of
\cite{mra-preunit}) where the action is $\varphi_{A\otimes V}=\mu_{A}\otimes V$.

From now on we consider quadruples $(A, V, \psi_{V}^{A}, \sigma_{V}^{A})$ where $A$
is an algebra, $V$ an object, $\psi_{V}^{A}$ satisfies (\ref{wmeas-wcp}) and
$\sigma_{V}^{A}:V\otimes V\rightarrow A\otimes V$ is a morphism in ${\mathcal C}$. For the idempotent morphism $\nabla_{A\otimes V}$ defined in (\ref{idem-wcp})
we denote by $A\times V$ the image of $\nabla_{A\otimes V}$, and
by $i_{A\otimes V}:A\times V\rightarrow A\otimes V$ and $p_{A\otimes
V}:A\otimes V\rightarrow A\times V$ the injection and the projection
associated to $\nabla_{A\otimes V}$.

We say that
$(A, V, \psi_{V}^{A},
\sigma_{V}^{A})$ satisfies the twisted condition if
\begin{equation}\label{twis-wcp}
(\mu_A\otimes V)\circ (A\otimes \psi_{V}^{A})\circ
(\sigma_{V}^{A}\otimes A) = (\mu_A\otimes V)\circ (A\otimes
\sigma_{V}^{A})\circ (\psi_{V}^{A}\otimes V)\circ (V\otimes
\psi_{V}^{A}).
\end{equation}
and   the  cocycle
condition holds if
\begin{equation}\label{cocy2-wcp}
(\mu_A\otimes V)\circ (A\otimes \sigma_{V}^{A}) \circ
(\sigma_{V}^{A}\otimes V) = (\mu_A\otimes V)\circ (A\otimes
\sigma_{V}^{A})\circ (\psi_{V}^{A}\otimes V)\circ
(V\otimes\sigma_{V}^{A}).
\end{equation}
By virtue of (\ref{twis-wcp}) and (\ref{cocy2-wcp}) we will consider from now
on, and without loss of generality, that
\begin{equation}
\label{idemp-sigma-inv}
\nabla_{A\otimes
V}\circ\sigma_{V}^{A} = \sigma_{V}^{A}
\end{equation}
holds for all quadruples $(A, V, \psi_{V}^{A}, \sigma_{V}^{A})$ {( see Proposition 3.7 of \cite{mra-preunit})}.

For $(A, V, \psi_{V}^{A}, \sigma_{V}^{A})$ define the product
\begin{equation}\label{prod-todo-wcp}
\mu_{A\otimes  V} = (\mu_A\otimes V)\circ (\mu_A\otimes
\sigma_{V}^{A})\circ (A\otimes \psi_{V}^{A}\otimes V)
\end{equation}
and let $\mu_{A\times V}$ be the product
\begin{equation}
\label{prod-wcp} \mu_{A\times V} = p_{A\otimes
V}\circ\mu_{A\otimes V}\circ (i_{A\otimes V}\otimes i_{A\otimes
V}).
\end{equation}
As a consequence of the twisted and the cocycle conditions, the product $\mu_{A\otimes V}$ is associative and normalized with respect to $\nabla_{A\otimes
V}$ (i.e. $\nabla_{A\otimes
V}\circ \mu_{A\otimes V}=\mu_{A\otimes V}=\mu_{A\otimes V}\circ (\nabla_{A\otimes
V}\otimes \nabla_{A\otimes
V}$)). Due to this normality condition, $\mu_{A\times V}$ is associative as well (Propostion 3.7 of \cite{mra-preunit}). Hence we define:
\begin{definicion}\label{wcp-def}
If $(A, V, \psi_{V}^{A}, \sigma_{V}^{A})$  satisfies (\ref{twis-wcp}) and (\ref{cocy2-wcp})
we say that $(A\otimes V, \mu_{A\otimes V})$ is a weak
crossed product.
\end{definicion}
Our next aim is to endow $A\times V$ with a unit, and hence with an algebra structure. As $A\times V$ is given as an image of an idempotent, it seems reasonable to use a preunit on $A\otimes V$ to obtain a unit on $A\times V$. In general, if $A$ is an algebra, $V$ an object in ${\mathcal C}$ and $m_{A\otimes V}$ is an associative product defined in
$A\otimes V$ a preunit $\nu:K\rightarrow A\otimes V$ is a morphism satisfying
\begin{equation}
m_{A\otimes V}\circ (A\otimes V\otimes \nu)=m_{A\otimes V}\circ (\nu\otimes
A\otimes V)=m_{A\otimes V}\circ (A\otimes V\otimes (m_{A\otimes V}\circ
(\nu\otimes \nu))).
\end{equation}
Associated to a preunit we obtain an idempotent morphism
$\nabla_{A\otimes
V}^{\nu}=m_{A\otimes V}\circ (A\otimes V\otimes \nu):A\otimes V\rightarrow
A\otimes V$. Take $A\times V$ the image of this idempotent, $p_{A\otimes V}^{\nu}$ the projection and $i_{A\otimes V}^{\nu}$ the injection. It is possible to endow $A\times V$ with an algebra structure whose product is $$m_{A\times V} = p_{A\otimes
V}^{\nu}\circ m_{A\otimes V}\circ (i_{A\otimes V}^{\nu}\otimes i_{A\otimes
V}^{\nu})$$
and whose unit is $\eta_{A\times V}=p_{A\otimes V}^{\nu}\circ \nu$ (see Proposition 2.5 of \cite{mra-preunit}). If moreover, $\mu_{A\otimes V}$ is left
$A$-linear for the actions $\varphi_{A\otimes V}=\mu_{A}\otimes V$, $\varphi_{A\otimes V\otimes
A\otimes V }=\varphi_{A\otimes V}\otimes  A\otimes V$ and normalized with respect to
$\nabla_{A\otimes
V}^{\nu}$,  the morphism
\begin{equation}
\label{beta-nu}
\beta_{\nu}:A\rightarrow A\otimes V,\; \beta_{\nu} =
(\mu_A\otimes V)\circ (A\otimes \nu)
\end{equation}
is
multiplicative and left $A$-linear for $\varphi_{A}=\mu_{A}$.

Although $\beta_{\nu}$ is not an algebra morphism, because
$A\otimes V$ is not an algebra, we have that $\beta_{\nu}\circ
\eta_A = \nu$, and thus the morphism $\bar{\beta_{\nu}} = p_{A\otimes
V}^{\nu}\circ\beta_{\nu}:A\rightarrow
A\times V$ is an algebra morphism.

In light of the considerations made in the last paragraph, and using the twisted and the cocycle conditions, we characterize weak crossed products with a preunit, and moreover we obtain an algebra structure on $A\times V$ (see \cite{mra-preunit}).

\begin{teorema}
\label{thm1-wcp} Let $A$ be an algebra, $V$ an object and
$m_{A\otimes V}:A\otimes V\otimes A\otimes V\rightarrow A\otimes V$ a
morphism of left $A$-modules  for the actions $\varphi_{A\otimes V}=\mu_{A}\otimes
V$, $\varphi_{A\otimes V\otimes A\otimes V }=\varphi_{A\otimes V}\otimes  A\otimes V$.

Then the following statements are equivalent:
\begin{itemize}
\item[(i)] The product $m_{A\otimes V}$ is associative with preunit
$\nu$ and normalized with respect to $\nabla_{A\otimes V}^{\nu}.$

\item[(ii)] There exist morphisms $\psi_{V}^{A}:V\otimes A\rightarrow A\otimes
V$, $\sigma_{V}^{A}:V\otimes V\rightarrow A\otimes V$ and $\nu:k\rightarrow
A\otimes V$ such that if $\mu_{A\otimes V}$ is the product defined
in (\ref{prod-todo-wcp}), the pair $(A\otimes V, \mu_{A\otimes
V})$ is a weak crossed product with $m_{A\otimes V} =
\mu_{A\otimes V}$ satisfying:
    \begin{equation}\label{pre1-wcp}
    (\mu_A\otimes V)\circ (A\otimes \sigma_{V}^{A})\circ
    (\psi_{V}^{A}\otimes V)\circ (V\otimes \nu) =
    \nabla_{A\otimes V}\circ
    (\eta_A\otimes V),
    \end{equation}
    \begin{equation}\label{pre2-wcp}
    (\mu_A\otimes V)\circ (A\otimes \sigma_{V}^{A})\circ
    (\nu\otimes V) = \nabla_{A\otimes V}\circ (\eta_A\otimes V),
    \end{equation}
    \begin{equation}\label{pre3-wcp}
(\mu_A\otimes V)\circ (A\otimes \psi_{V}^{A})\circ (\nu\otimes A)
= \beta_{\nu},
\end{equation}
\end{itemize}
where $\beta_{\nu}$ is the morphism defined in
(\ref{beta-nu}). In this case $\nu$ is a preunit for $\mu_{A\otimes
V}$, the idempotent morphism of the weak crossed product
$\nabla_{A\otimes V}$ is the idempotent $\nabla_{A\otimes
V}^{\nu}$, and we say that the pair $(A\otimes V, \mu_{A\otimes V})$ is a
weak crossed product with preunit $\nu$.
\end{teorema}
As a corollary of  Theorem \ref{thm1-wcp}  we obtain:
\begin{corol}\label{corol-wcp}
If $(A\otimes V, \mu_{A\otimes V})$ is a weak crossed product with
preunit $\nu$, then $A\times V$ is an algebra with the product
defined in (\ref{prod-wcp}) and unit $\eta_{A\times V}=p_{A\otimes
V}\circ\nu$.
\end{corol}

In the rest of the section we obtain some new results about weak crossed products. The following theorem gives necessary and sufficient conditions for an algebra $B$ to be isomorphic to the algebra derived from a weak crossed product with preunit.

\begin{teorema}
\label{uni-1}
Let $B$ be an algebra in ${\mathcal C}$. Then the following are equivalent:

\begin{itemize}

\item[(i)] There exist a weak crossed product $(A\otimes V, \mu_{A\otimes V})$ with
preunit $\nu$ and an isomorphism of algebras $\omega:A\times V\rightarrow B$.
\item[(ii)] There exist an algebra $A$, an object $V$, morphisms $$i_{A}:A\rightarrow
B,\;\;\;\; i_{V}:V\rightarrow B,\;\;\;\;\nabla_{A\otimes V}: A\otimes V\rightarrow A\otimes
V, \;\;\;\;\omega:A\times V\rightarrow
B$$
such that $i_{A}$ is an algebra morphism, $\nabla_{A\otimes V}$ is an idempotent
morphism  of left $A$-modules for the action $\varphi_{A\otimes V}=\mu_{A}\otimes V$ and $\omega$ is an isomorphism such that
$$\omega\circ p_{A\otimes V}=\mu_{B}\circ (i_{A}\otimes i_{V})$$
 where $A\times V$ is the image of
$\nabla_{A\otimes V}$ and $p_{A\otimes V}$ is the associated projection.

\item[(iii)] There exist an algebra $A$, an object $V$ and morphisms $i_A:A\rightarrow B$, $i_V:V\rightarrow B$ and $\hat{\omega}:B\rightarrow A\otimes V$ that satisfy:
\begin{itemize}
\item[(iii-1)] $i_A$ is a morphism of algebras.
\item[(iii-2)] $\hat{\omega}$ is a morphism of left $A$-modules for $\varphi_{B} = \mu_B\circ (i_A\otimes B)$ and $\varphi_{A\otimes V} = \mu_A\otimes V$.
\item[(iii-3)] $\mu_B\circ (i_A\otimes i_V)\circ \hat{\omega} = id_B$.
\end{itemize}

\end{itemize}

\end{teorema}

\begin{dem}
First we prove (i) $\Rightarrow$ (ii). Suppose that  $(A\otimes V, \mu_{A\otimes
V})$ is a weak crossed product with preunit $\nu:K\rightarrow A\otimes V$. Define $i_{A}:A\rightarrow B$ by
\begin{equation}
\label{iA}
i_{A}=\omega\circ p_{A\otimes V}\circ \beta_{\nu}:A\rightarrow B
\end{equation}
where $p_{A\otimes V}$ is the projection associated to the idempotent and
$\beta_{\nu}$ the morphism defined in (\ref{beta-nu}). As in this case $p_{A\otimes V} = p_{A\otimes V}^{\nu}$, then $p_{A\otimes V}\circ \beta_{\nu} = \bar{\beta}_{\nu}$ is a morphism of algebras, and thus $i_A$ is a morphism of algebras.

Define  $i_{V}: V\rightarrow B$ by
\begin{equation}
\label{iV}
i_{V}=\omega\circ p_{A\otimes V}\circ (\eta_{A}\otimes V).
\end{equation}
Then, $\omega\circ p_{A\otimes V}=\mu_{B}\circ (i_{A}\otimes i_{V})$ holds if and only if
$p_{A\times V}=\omega^{-1}\circ \mu_{B}\circ (i_{A}\otimes i_{V})$. This last equality follows by:

\begin{itemize}
\item[ ]$\hspace{0.38cm} \omega^{-1}\circ  \mu_{B}\circ (i_{A}\otimes i_{V})  $

\item[ ]$= p_{A\otimes V} \circ \mu_{A\otimes V}\circ ((\nabla_{A\otimes V}\circ
\beta_{\nu})\otimes (\nabla_{A\otimes V}\circ (\eta_{A}\otimes V))$

\item[ ]$= p_{A\otimes V} \circ \mu_{A\otimes V}\circ ( \beta_{\nu}\otimes
\eta_{A}\otimes V)$

\item[ ]$=p_{A\otimes V} \circ (\mu_{A}\otimes V)\circ (A\otimes
\sigma_{V}^{A})\circ ((\nabla_{A\otimes V}\circ \beta_{\nu})\otimes V)$

\item[ ]$=p_{A\otimes V} \circ (\mu_{A}\otimes V)\circ (A\otimes
\sigma_{V}^{A})\circ (\beta_{\nu}\otimes V)$

\item[ ]$=p_{A\otimes V} \circ (\mu_{A}\otimes V)\circ (A\otimes ((\mu_{A}\otimes
V)\circ (A\otimes \sigma_{V}^{A})\circ (\nu\otimes V)))$

\item[ ]$=p_{A\otimes V} \circ (\mu_{A}\otimes V)\circ (A\otimes (\nabla_{A\otimes
V}\circ (\eta_{A}\otimes V)))$

\item[ ]$= p_{A\otimes V} \circ \nabla_{A\otimes V}$

\item[ ]$=p_{A\otimes V}.  $

\end{itemize}

The first equality is consequence of the condition of algebra morphism for $\omega$ and the
definition of $\mu_{A\times V}$, the second one of the normalized character of
$\mu_{A\otimes V}$, the third one of (\ref{prod-todo-wcp}), the fourth one of
(\ref{pre3-wcp}), the fifth one of the associativity of $\mu_{A}$, the sixth one of
(\ref{pre2-wcp}), the seventh one of the left $A$-module condition for the idempotent
and finally the eighth one of the properties of $i_{A\otimes V}$ and $p_{A\otimes
V}$.

The proof for (ii)$\Rightarrow$(i) is the following:
Suppose now that the conditions of the second statement are fulfilled, and define the product
\begin{equation}
\label{apoyo1}
\mu_{A\times V}=\omega^{-1}\circ \mu_{B}\circ (\omega\otimes \omega).
\end{equation}
that induces an algebra structure on $A\times V$ with unit $\eta_{A\times V}= \omega^{-1}\circ \eta_B$, and forces $\omega$ to be an isomorphism of algebras. Moreover we can define on $A\otimes V$ the multiplication given by
\begin{equation}
\label{apoyo}
\mu_{A\otimes V}=i_{A\otimes V}\circ \mu_{A\times V}\circ ( p_{A\otimes V}\otimes
p_{A\otimes V}).
\end{equation}
Now to complete the proof, it is enough to check that $(A\otimes V, \mu_{A\otimes V})$ is under the assumptions of (i) of Theorem \ref{thm1-wcp}, so we can assure that $(A\otimes V, \mu_{A\otimes V})$ is induced by a weak crossed product with preunit.
First observe that $\mu_{A\otimes V}$ is associative for being (\ref{apoyo1}) associative, and moreover, as a consequence of the properties of $i_{A\otimes V}$ and $p_{A\otimes V}$ the product is normalized with respect to $\nabla_{A\otimes V}$. To define a preunit on $(A\otimes V, \mu_{A\otimes V})$ consider
$$\nu=i_{A\otimes V}\circ \eta_{A\times V}:K\rightarrow A\otimes V.$$
This morphism is a preunit for $(A\otimes V, \mu_{A\otimes V})$ as, on the one hand $\mu_{A\otimes V}\circ (A\otimes V\otimes \nu)=\mu_{A\otimes
V}\circ (\nu\otimes A\otimes V)$
and on the other hand $\mu_{A\otimes V}\circ (\nu\otimes\nu) = \nu$.
It is easy to check that
\begin{equation}
\label{idempot-preunit}
\nabla_{A\otimes V} = \mu_{A\otimes V}\circ (A\otimes V\otimes \nu).
\end{equation}
 Hence $\nabla_{A\otimes V}^{\nu}=\nabla_{A\otimes V}$, and thus $\mu_{A\otimes V}$ is normalized with respect to $\nabla_{A\otimes V}^{\nu}$. To finish the proof it just remain to check that $\mu_{A\otimes V}$ is a morphism of left $A$-modules. In order to obtain this result note that, as $i_A$ is an algebra morphism:
\begin{equation}
\label{apoyo2}
\mu_{A\otimes V}\circ ((i_{A\otimes V}\circ \omega^{-1}\circ i_{A})\otimes A\otimes
V)=\nabla_{A\otimes V}\circ (\mu_{A}\otimes V).
\end{equation}

Then,  if $\varphi_{A\otimes V}$ and $\varphi_{A\otimes V\otimes A\otimes V}$ are
the actions used in Theorem \ref{thm1-wcp}, we have

\begin{itemize}
\item[ ]$\hspace{0.38cm} \mu_{A\otimes V}\circ \varphi_{A\otimes V\otimes A\otimes
V}  $

\item[ ]$= i_{A\otimes V}\circ \omega^{-1}\circ \mu_{B}\circ ((\mu_{B}\circ
(i_{A}\otimes i_{V})\circ  \nabla_{A\otimes V}\circ (\mu_{A}\otimes V))\otimes
(\mu_{B}\circ (i_{A}\otimes i_{V})\circ  \nabla_{A\otimes V}))$

\item[ ]$=  i_{A\otimes V}\circ \omega^{-1}\circ \mu_{B}\circ ((\mu_{B}\circ
(i_{A}\otimes i_{V})\circ (\mu_{A}\otimes V)\circ (A\otimes \nabla_{A\otimes
V}))\otimes (\mu_{B}\circ (i_{A}\otimes i_{V})\circ  \nabla_{A\otimes V})) $

\item[ ]$= i_{A\otimes V}\circ \omega^{-1}\circ \mu_{B}\circ (i_{A}\otimes
(\mu_{B}\circ (\omega\otimes \omega)\circ (p_{A\otimes V} \otimes p_{A\otimes V})))$

\item[ ]$= i_{A\otimes V}\circ \omega^{-1}\circ \mu_{B}\circ ((\omega\circ
\omega^{-1}\circ i_{A})\otimes (\omega\circ  \mu_{A\times V}\circ (p_{A\otimes V}
\otimes p_{A\otimes V})))$

\item[ ]$=i_{A\otimes V}\circ \mu_{A\times V} \circ ((\omega^{-1}\circ i_{A})\otimes
( \mu_{A\times V}\circ (p_{A\otimes V} \otimes p_{A\otimes V})))  $

\item[ ]$= \mu_{A\otimes V} \circ ((i_{A\otimes V}\circ \omega^{-1}\circ
i_{A})\otimes \mu_{A\otimes V})  $

\item[ ]$= \nabla_{A\otimes V}\circ (\mu_{A}\otimes V)\circ (A\otimes \mu_{A\otimes
V}) $

\item[ ]$=(\mu_{A}\otimes V)\circ (A\otimes ( \nabla_{A\otimes V}\circ \mu_{A\otimes
V}))   $

\item[ ]$=(\mu_{A}\otimes V)\circ (A\otimes \mu_{A\otimes V})$

\item[ ]$=\varphi_{A\otimes V}\circ (A\otimes \mu_{A\otimes V}).$

\end{itemize}

The first equality follows by definition, the second one by the left $A$-module
condition for $\nabla_{A\otimes V}$, the third one by the associativity of $\mu_{B}$
and  the condition of algebra morphism of $i_{A}$, the fourth and the fifth ones  by
(\ref{apoyo1}), the sixth one by (\ref{apoyo}), the seventh one by (\ref{apoyo2}),
the eighth one by  the left $A$-module condition for $\nabla_{A\otimes V}$, the
ninth one by the normal condition and finally, the tenth one by definition.

Therefore, Theorem \ref{thm1-wcp} assures that $(A\otimes V,\mu_{A\otimes V})$ is a
weak crossed with preunit $\nu$ and this finishes the proof of (ii)$\Rightarrow$(i).

To prove (i)$\Rightarrow$(iii) define $i_{A}$ and $i_{V}$ as in (i)$\Rightarrow$(ii) and put  $\hat{\omega} = i_{A\otimes V}\circ \omega^{-1}:B\rightarrow A\otimes V$. Then,
$$\mu_B\circ (i_A\otimes i_V)\circ \hat{\omega} = \omega\circ p_{A\otimes V}\circ i_{A\otimes V} \circ \omega^{-1} = id_B.$$
 and $\hat{\omega}$ is a morphism of left $A$-modules because:

\begin{itemize}

\item[ ]$\hspace{0.38cm} \hat{\omega}\circ \varphi_{B}$

\item[ ]$= i_{A\otimes V}\circ \omega^{-1}\circ \mu_B\circ (i_A\otimes B) $

\item[ ]$= i_{A\otimes V}\circ \mu_{A\times V}\circ ((\omega^{-1}\circ i_A)\otimes \omega^{-1})$

\item[ ]$= i_{A\otimes V}\circ \mu_{A\times V}\circ ((p_{A\otimes V}\circ \beta_{\nu})\otimes \omega^{-1})$

\item[ ]$= \mu_{A\otimes V}\circ (\beta_{\nu}\otimes (i_{A\otimes V}\circ \omega^{-1}))$

\item[ ]$= (\mu_A\otimes V)\circ (A\otimes (\mu_{A\otimes V}\circ (\nu\otimes (i_{A\otimes V}\circ \omega^{-1}))))$

\item[ ]$= \varphi_{A\otimes V}\circ (A\otimes \hat{\omega})$

\end{itemize}

\noindent where the first equality follows by the definition of $\hat{\omega}$, the second one is consequence of being $\omega^{-1}$ of algebras, the third one follows by the definition of $i_A$, the fourth one by the definition of $\mu_{A\times V}$, by the normality of $\mu_{A\otimes V}$,  and by $\nabla_{A\otimes V}\circ \beta_{\nu} = \beta_{\nu}$, the fifth one is consequence of being $\mu_{A\otimes V}$ of left $A$-modules and the last one is consequence of being $\nabla_{A\otimes V} = \nabla_{A\otimes V}^{\nu}$ and $\nabla_{A\otimes V}\circ i_{A\otimes V}= i_{A\otimes V}$.

Finally, we prove (iii)$\Rightarrow$(i). Now we have a morphism of algebras $i_A:A\rightarrow B$ and a morphism $i_V:V\rightarrow B$ such that there exists $\hat{\omega}:B\rightarrow A\otimes V$ of left $A$-modules that satisfies $\mu_B\circ (i_A\otimes i_V)\circ \hat{\omega} = id_B$. It is clear that
$$\Omega = \hat{\omega}\circ \mu_B\circ (i_A\otimes i_V):A\otimes V\rightarrow A\otimes V$$
is an idempotent morphism  such that
\begin{itemize}

\item[ ]$\hspace{0.38cm} \varphi_{A\otimes V}\circ (A\otimes \Omega)$

\item[ ]$=  \hat{\omega} \circ \mu_{B}\circ (i_{A}\otimes (\mu_{B}\circ (i_{A}\otimes i_{V})))$

\item[ ]$= \hat{\omega} \circ \mu_{B}\circ ((\mu_{B}\circ (i_{A}\otimes i_{A}))\otimes i_{V})$

\item[ ]$= \hat{\omega} \circ \mu_{B}\circ ((i_{A}\circ \mu_{A})\otimes i_{V})$

\item[ ]$= \Omega\circ \varphi_{A\otimes V}$

\end{itemize}
 and then is a morphism of left $A$-modules.

Define the product $\mu_{A\otimes V}:A\otimes V\otimes A\otimes V\rightarrow A\otimes V$ by
$$\mu_{A\otimes V}=\hat{\omega}\circ \mu_{B}\circ ( \bar{\omega}\otimes  \bar{\omega})$$
where $\bar{\omega} = \mu_B\circ (i_A\otimes i_V)$. As $\bar{\omega}\circ \hat{\omega} = id_B$, and as $i_A$ is of algebras, this product is associative, of left $A$-modules and with preunit $\nu = \hat{\omega}\circ \eta_B$. It is also easy to check that $\Omega = \nabla_{A\otimes V}^{\nu}$ and $\Omega\circ \mu_{A\otimes V} = \mu_{A\otimes V}= \mu_{A\otimes V}\circ (\Omega\otimes \Omega)$, that is, $\mu_{A\otimes V}$ is normalized. Hence we are under the conditions of Theorem \ref{thm1-wcp}, and $(A\otimes V, \mu_{A\otimes V})$ is a weak crossed product with preunit $\nu$, associated idempotent morphism $\nabla_{A\otimes V} = \Omega$,
$A\times V$ the image of the idempotent and $p_{A\otimes V}$ and $i_{A\otimes V}$
the projection and the injection respectively.

Put $\omega= \bar{\omega}\circ i_{A\otimes V}:A\times V\rightarrow B$. Then, $\omega$ is an isomorphism with inverse $\omega^{-1}=
p_{A\otimes V}\circ \hat{\omega}$ and the following diagram is commutative

\[\xymatrix{ & A\times V\ar @/^/[dd]^(0.35){\omega}\ar[dr]^-{i_{A\otimes V}} & \\
A\otimes V\ar[rr]_(0.45){\Omega}\ar[ur]^-{p_{A\otimes V}}\ar[dr]_-{\bar{\omega}} && A\otimes V\\
& B\ar[ur]_-{\hspace{0.2cm}\hat{\omega}}}\]

Moreover, $\omega $ is an algebra morphism because
$$\omega\circ \eta_{A\times V}=\mu_B\circ (i_A\otimes i_V)\circ \nabla_{A\otimes V}\circ \nu=\mu_B\circ (i_A\otimes i_V)\circ \hat{\omega}\circ\eta_{B }=\eta_{B}$$
and
$$\omega^{-1}\circ\mu_{B}\circ (\omega\otimes \omega)$$
$$= p_{A\otimes V}\circ \hat{\omega}\circ \mu_{B}\circ ( \bar{\omega}\otimes  \bar{\omega})\circ (i_{A\otimes V}\otimes i_{A\otimes V})= p_{A\otimes V}\circ \mu_{A\otimes V}\circ (i_{A\otimes V}\otimes i_{A\otimes V})
$$
$$=\mu_{A\times V}.$$
\end{dem}

\begin{nota}
\label{proyeccion}
Note that, if the conditions (ii) or (iii) of Theorem \ref{uni-1} hold it is possible
to characterize the isomorphism $\omega$ as
$$\omega=\mu_{B}\circ (i_{A}\otimes i_{V})\circ i_{A\otimes V}.$$

If we are under the conditions (iii) of Theorem \ref{uni-1}, we know that
$A\otimes V$ has a weak crossed product structure with preunit $\nu =
\hat{\omega}\circ \eta_B$. Observe that this fact implies that

\begin{equation}
\label{unidad-uni}
\bar{\omega}\circ \nu = \mu_B\circ (i_A\otimes i_V)\circ
\nu = \eta_B.
\end{equation}

If $\omega:A\times V\rightarrow B$ is the isomorphism of algebras, in light of the equality  $\omega^{-1}=
p_{A\otimes V}\circ \hat{\omega}$ and the fact that $i_A$ is a morphism of algebras we
obtain that
\[\omega\circ p_{A\otimes V}\circ \beta_{\nu} = i_A\]
and
\[\omega\circ p_{A\otimes V}\circ (\eta_A\otimes V) = i_V.\]
Thus the morphisms $i_A$ and $i_V$ have always this explicit expression.
\end{nota}

\begin{nota} \label{universal-1}

Suppose that $B$ satisfies the conditions of Theorem \ref{uni-1}, and let $i_{A}$,
$i_{V}$ be the morphisms defined in (\ref{iA}) and (\ref{iV}). By Theorem \ref{thm1-wcp} of \cite{mra-preunit} we know that
\begin{equation}\label{fi-wcp}
\psi_{V}^{A} = \mu_{A\otimes V}\circ (\eta_A\otimes V\otimes
\beta_{\nu})
\end{equation}
\begin{equation}\label{sigma-wcp}
\sigma_{V}^{A} = \mu_{A\otimes V}\circ (\eta_A\otimes V\otimes
\eta_A\otimes V)
\end{equation}
and then, by the condition of algebra morphism for $\omega$ and the equality
$\nabla_{A\otimes V}\circ \psi_{V}^{A} =  \psi_{V}^{A}$ we obtain
\begin{equation}
\label{psi-uni}
\mu_{B}\circ (i_{A}\otimes i_{V})\circ \psi_{V}^{A}=\mu_{B}\circ (i_{V}\otimes i_{A})
\end{equation}
and
\begin{equation}
\label{sigma-uni}
\mu_{B}\circ (i_{A}\otimes i_{V})\circ \sigma_{V}^{A}=\mu_{B}\circ (i_{V}\otimes
i_{V}).
\end{equation}

The conditions (\ref{unidad-uni}), (\ref{psi-uni}) and (\ref{sigma-uni}) are closely
related to the characterization of weak crossed products as universal
constructions. In the following result we obtain this universal property explicitly.

\end{nota}

\begin{teorema}
\label{universal}
Let $(A\otimes V, \mu_{A\otimes V})$ be a weak crossed product  with preunit $\nu$,
and let $B$ be an algebra. Suppose that there exist morphisms $i_{A}:A\rightarrow B$
and $i_{V}:V\rightarrow B$ such that $i_{A}$ is an algebra morphism. Then the following assertions are equivalent:

\begin{itemize}

\item[(i)] There exists an unique algebra morphism $\omega:A\times V\rightarrow B$ that makes the following diagram commutative:

\begin{center}
\scalebox{0.85}
{
\begin{pspicture}(0,-3.2829688)(8.502812,3.2829688)
\psline[linewidth=0.04cm,arrowsize=0.05291667cm
2.0,arrowlength=1.4,arrowinset=0.4]{->}(0.7009375,-0.27546874)(4.4209375,-0.25546876)
\psline[linewidth=0.04cm,arrowsize=0.05291667cm
2.0,arrowlength=1.4,arrowinset=0.4]{<-}(0.2209375,0.10453125)(6.0409374,2.7645311)
\psline[linewidth=0.04cm,arrowsize=0.05291667cm
2.0,arrowlength=1.4,arrowinset=0.4]{<-}(4.9609375,0.06453125)(6.4209375,2.5845313)
\psline[linewidth=0.04cm,arrowsize=0.05291667cm
2.0,arrowlength=1.4,arrowinset=0.4]{->}(5.9809375,-3.2354689)(0.2409375,-0.67546874)
\psline[linewidth=0.04cm,arrowsize=0.05291667cm
2.0,arrowlength=1.4,arrowinset=0.4]{->}(6.2209377,-2.8554688)(4.9409375,-0.53546876)
\usefont{T1}{ptm}{m}{n}
\rput(-0.25,-0.25){$A\times V$}
\usefont{T1}{ptm}{m}{n}
\rput(4.9,-0.2){$B$}
\usefont{T1}{ptm}{m}{n}
\rput(6.402344,2.9){$A$}
\usefont{T1}{ptm}{m}{n}
\rput(6.362344,-3.2054688){$V$}
\usefont{T1}{ptm}{m}{n}
\rput(6.25,1.3745313){$i_{A}$}
\usefont{T1}{ptm}{m}{n}
\rput(6.1,-1.4654688){$i_{V}$}
\usefont{T1}{ptm}{m}{n}
\rput(2.5623438,0.03453125){$\omega$}
\usefont{T1}{ptm}{m}{n}
\rput(2.5,2.0945313){$p_{A\otimes V}\circ \beta_{\nu}$}
\usefont{T1}{ptm}{m}{n}
\rput(1.3,-2.2254686){$p_{A\otimes V}\circ (\eta_{A}\otimes V)$}
\end{pspicture}
}
\end{center}

\item[(ii)] The
equalities  (\ref{unidad-uni}), (\ref{psi-uni}) and (\ref{sigma-uni}) hold.

\end{itemize}

\end{teorema}

\begin{dem}
First we prove (i)$\Rightarrow$(ii). Suppose that $\omega:A\times V\rightarrow B$ is a morphism of algebras commuting the diagram. Then
\begin{itemize}

\item[ ]$\hspace{0.38cm} \mu_B\circ (i_{A}\otimes i_{V})$

\item[ ]$= \mu_{B}\circ ((\omega\circ p_{A\otimes V}\circ \beta_{\nu})\otimes (\omega\circ p_{A\otimes V}\circ (\eta_{A}\otimes V)))  $

\item[ ]$= \omega\circ \mu_{A\times V}\circ (( p_{A\otimes V}\circ \beta_{\nu})\otimes ( p_{A\otimes V}\circ (\eta_{A}\otimes V)))$

\item[ ]$=\omega\circ p_{A\otimes V}\circ \mu_{A\otimes V}\circ ( \beta_{\nu}\otimes \eta_{A}\otimes V) $

\item[ ]$= \omega\circ p_{A\otimes V}\circ \mu_{A\otimes V}\circ (\mu_{A}\otimes V\otimes A\otimes V)\circ (A\otimes \nu\otimes \eta_{A}\otimes V)$

\item[ ]$= \omega\circ p_{A\otimes V}\circ (\mu_{A}\otimes V)\circ (A\otimes (\mu_{A\otimes V}\circ (\nu\otimes \eta_{A}\otimes V)))  $

\item[ ]$= \omega\circ p_{A\otimes V}\circ (\mu_{A}\otimes V)\circ (A\otimes (\nabla_{A\otimes V}\circ ( \eta_{A}\otimes V))) $

\item[ ]$=\omega\circ p_{A\otimes V}$

\end{itemize}

where the first equality follows by the commutativity of the diagram, the second one by the condition of algebra morphism for $\omega$, the third one by the definition of $\mu_{A\times V}$, the fourth one is a consequence  of the definition of $\beta_{\nu}$, the fifth one uses that $\mu_{A\otimes V}$ is a left $A$-module morphism, the sixth one follows by $\nabla_{A\otimes V}=\nabla_{A\otimes V}^{\nu}$ and finally, the seventh one uses that $\nabla_{A\otimes V}$ is a left $A$-module morphism.

Then, as a consequence,  we have the following:
$$ \mu_B\circ (i_{A}\otimes i_{V})\circ \nu=\omega\circ p_{A\otimes V}\circ \nu=\omega\circ \eta_{A\otimes V}=\eta_{B}$$
and (\ref{unidad-uni}) holds.
On the other hand,
 $$\mu_B\circ (i_{A}\otimes i_{V})\circ \psi_{V}^{A}=\omega\circ p_{A\otimes V}\circ \mu_{A\otimes V}\circ (\eta_{A}\otimes V\otimes \beta_{\nu})=\omega\circ \mu_{A\times V}\circ ((p_{A\otimes V}\circ (\eta_{A}\otimes V))\otimes (p_{A\otimes V}\circ \beta_{\nu}))$$
$$= \mu_{B}\circ ((\omega\circ p_{A\otimes V}\circ (\eta_{A}\otimes V))\otimes (\omega\circ p_{A\otimes V}\circ \beta_{\nu}))=\mu_{B}\circ (i_{V}\otimes i_{A})$$
and (\ref{psi-uni}) holds.

Finally, we obtain (\ref{sigma-uni}) by  similar computations.
$$\mu_B\circ (i_{A}\otimes i_{V})\circ \sigma_{V}^{A}=\omega\circ p_{A\otimes V}\circ \mu_{A\otimes V}\circ (\eta_{A}\otimes V\otimes \eta_{A}\otimes V)$$
$$=\omega\circ \mu_{A\times V}\circ ((p_{A\otimes V}\circ (\eta_{A}\otimes V))\otimes (p_{A\otimes V}\circ (\eta_{A}\otimes V)))$$
$$= \mu_{B}\circ ((\omega\circ p_{A\otimes V}\circ (\eta_{A}\otimes V))\otimes (\omega\circ p_{A\otimes V}\circ (\eta_{A}\otimes V)))=\mu_{B}\circ (i_{V}\otimes i_{V})$$

The proof for (ii)$\Rightarrow$(i) is the following. Define $\omega=\mu_{B}\circ (i_{A}\otimes i_{V})\circ i_{A\otimes V}:A\times
V\rightarrow B$. First we need to prove that this morphism is multiplicative, and this fact follows by:
\begin{itemize}
\item[ ]$\hspace{0.38cm}\omega\circ \mu_{A\times V}   $

\item[ ]$=\mu_{B}\circ ((\mu_{B}\circ (i_{A}\otimes i_{A}))\otimes (\mu_{B}\circ
(i_{A}\otimes i_{V})\circ \sigma_{V}^{A}))\circ (A\otimes \psi_{V}^{A}\otimes
V)\circ (i_{A\otimes V}\otimes i_{A\otimes V}) $

\item[ ]$= \mu_{B}\circ (B\otimes \mu_{B})\circ (i_{A}\otimes (\mu_{B}\circ
(i_{A}\otimes i_{V})\circ \psi_{V}^{A})\otimes i_{V})\circ (i_{A\otimes V}\otimes
i_{A\otimes V}) $

\item[ ]$=\mu_{B}\circ (\omega\otimes \omega).  $
\end{itemize}

The first equality follows by the associativity of $\mu_{B}$ and the condition of
algebra morphism for $i_{A}$, the second one by (\ref{sigma-uni}) and the
associativity of $\mu_{B}$. The third one is a consequence of (\ref{psi-uni}) as
well as the associativity of $\mu_{B}$.

Now by (\ref{unidad-uni}) we obtain that $\omega\circ \eta_{A\otimes V}=\eta_{B}$ and therefore $\omega$ is an algebra morphism.

Moreover, equalities (\ref{unidad-uni}), (\ref{psi-uni}) and the
condition of algebra morphism for $i_{A}$ forces the diagram to be commutative. Indeed:

\begin{itemize}
\item[ ]$\hspace{0.38cm}\omega\circ p_{A\otimes V} \circ \beta_{\nu}  $

\item[ ]$=\mu_{B}\circ ((i_{A}\circ \mu_{A})\otimes i_{V})\circ (A\otimes
 \nu)  $

\item[ ]$= \mu_{B}\circ (i_{A}\otimes ( \mu_{B}\circ (i_{A}\otimes i_{V})\circ
 \nu))$

\item[ ]$=\mu_{B}\circ (i_{A}\otimes \eta_{B})$

\item[ ]$=i_{A}$
\end{itemize}

and

\begin{itemize}
\item[ ]$\hspace{0.38cm}\omega\circ p_{A\otimes V} \circ (\eta_{A}\otimes V) $

\item[ ]$=\mu_{B}\circ (i_{A}\otimes i_{V})\circ \nabla_{A\otimes V}\circ
(\eta_{A}\otimes V)  $

\item[ ]$= \mu_{B}\circ (i_{A}\otimes i_{V})\circ \psi_{V}^{A}\circ (V\otimes
\eta_{A})$

\item[ ]$=\mu_{B}\circ (i_{V}\otimes (i_{A}\circ \eta_{A}))$

\item[ ]$=\mu_{B}\circ (i_{V}\otimes \eta_{B})$

\item[ ]$=i_{V}.$
\end{itemize}

Finally, suppose that $\theta:A\times V\rightarrow B$ is a morphism that satisfies the same conditions as $\omega$. Then:

\begin{itemize}
\item[ ]$\hspace{0.38cm}\omega$

\item[ ]$=\mu_{B}\circ (i_{A}\otimes i_{V})\circ i_{A\otimes V} $

\item[ ]$=\mu_{B}\circ ( (\theta\circ p_{A\otimes V} \circ \beta_{\nu})\otimes
(\theta\circ p_{A\otimes V} \circ (\eta_{A}\otimes V)) \circ i_{A\otimes V}  $

\item[ ]$=\theta \circ p_{A\otimes V}\circ \mu_{A\otimes V}\circ (\beta_{\nu}\otimes
\eta_{A}\otimes V)\circ i_{A\otimes V}  $

\item[ ]$=\theta \circ p_{A\otimes V}\circ \nabla_{A\otimes V}\circ i_{A\otimes V} $

\item[ ]$=\theta$
\end{itemize}

and therefore $\omega$ is unique.
\end{dem}

\begin{nota} Theorems \ref{uni-1} and \ref{universal} are the generalization to the
weak case of Propositions 2.3 and 2.4 of \cite{bes-drab2}. In these results Bespalov
and Drabant described equivalent projection and injection conditions for an algebra to be isomorphic to a
cross product algebra and characterized the universal construction in
this setting.
\end{nota}

\begin{nota}\label{VA}
In the previous results we develop a theory of crossed products in $A\otimes
 V$. In a symmetric  way it is possible to obtain similar results
for $V\otimes A$. In this new case we must work with morphisms
$\psi_{A}^{V}:A\otimes V\rightarrow V\otimes A$,
$\sigma_{A}^{V}:V\otimes V\rightarrow V\otimes A$ and the associated
idempotent, $\nabla_{V\otimes A}:V\otimes A\rightarrow V\otimes A$,
is defined by
\begin{equation}
\label{otra-nabla}
\nabla_{V\otimes A}=(V\otimes \mu_{A})\circ ((\psi_{A}^{V}\circ (\eta_{A}\otimes
V))\otimes A).
\end{equation}
The induced product $\mu_{V\otimes A}:V\otimes A\otimes V\otimes A\rightarrow
V\otimes A$ is
\begin{equation}
\label{otro-prod}
\mu_{V\otimes A}=(V\otimes \mu_{A})\circ (\sigma_{A}^{V}\otimes \mu_{A})\circ
(V\otimes \psi_{A}^{V}\otimes A)
\end{equation}
and the preunit is a morphism $\nu:K\rightarrow V\otimes A$.
\end{nota}

The next proposition shows that under suitable conditions it is possible to find a weak crossed product over $V\otimes A$ once we have one defined on $A\otimes
V$.

\begin{prop}
\label{pro-copro}
Let $(A\otimes V,\mu_{A\otimes V})$ be a weak crossed product with associated
morphisms $\psi_{V}^{A}:V\otimes A\rightarrow A\otimes V$, $\sigma_{V}^{A}:V\otimes
V\rightarrow A\otimes V$. Let
$(A,V, \psi_{A}^{V}:A\otimes V\rightarrow V\otimes A)$ be a triple satisfying
\begin{equation}
\label{idem-sim}
(V\otimes \mu_{A})\circ (\psi_{A}^{V}\otimes A)\circ (A\otimes
\psi_{A}^{V})=\psi_{A}^{V}\circ (\mu_{A}\otimes V).
\end{equation}
Let $\nabla_{A\otimes V}$ be the idempotent of the weak crossed product $(A\otimes V, \mu_{A\otimes V})$ and $\nabla_{V\otimes A}$ a morphism given as in (\ref{otra-nabla}). Then if
\begin{equation}
\label{idem-psiAV}
\psi_{V}^{A}\circ \psi_{A}^{V}=\nabla_{A\otimes V}
\end{equation}
 and
\begin{equation}
\label{idem-psiVA}
\psi_{A}^{V}\circ \psi_{V}^{A}=\nabla_{V\otimes A},
\end{equation}
hold, the pair $(V\otimes A,\mu_{V\otimes A})$ is a weak crossed product with associated
morphisms $\psi_{A}^{V}$, $\sigma_{A}^{V}=\psi_{A}^{V}\circ \sigma_{V}^{A}:V\otimes
V\rightarrow V\otimes A$ and the product $\mu_{V\otimes A}$ defined in
(\ref{otro-prod}) satisfies
\begin{equation}
\label{idem-pro}
\mu_{V\otimes A}=\psi_{A}^{V}\circ \mu_{A\otimes V}\circ (\psi_{V}^{A}\otimes
\psi_{V}^{A}).
\end{equation}

Moreover, if $(A\otimes V,\mu_{A\otimes V})$ is a weak crossed product with preunit
$\nu:K\rightarrow A\otimes V$, the pair $(V\otimes A,\mu_{V\otimes A})$ is a weak
crossed product with preunit $\upsilon=\psi_{A}^{V}\circ \nu:
K\rightarrow V\otimes A$.
\end{prop}

\begin{dem}
First note that $\nabla_{V\otimes A}\circ \sigma_{A}^{V}=\sigma_{A}^{V}$ because by
(\ref{idem-psiAV}) and (\ref{idem-psiVA}),
$$\nabla_{V\otimes A}\circ \sigma_{A}^{V}=\psi_{A}^{V}\circ \psi_{V}^{A}\circ
\psi_{A}^{V}\circ \sigma_{V}^{A}=\psi_{A}^{V}\circ \nabla_{A\otimes V}\circ
\sigma_{V}^{A}=\psi_{A}^{V}\circ \sigma_{V}^{A} =\sigma_{A}^{V}.$$

Before proving the twisted and the cocycle conditions for $\psi_A^V$ and $\sigma_A^V$ observe that by (\ref{idem-sim}) we obtain:
\begin{equation}
\label{AV-axu1}
\psi_{A}^{V}\circ (\mu_{A}\otimes V)\circ (A\otimes \psi_{V}^{A})\circ (\sigma_{V}^{A}\otimes A)=(V\otimes
\mu_{A})\circ (\sigma_{A}^{V}\otimes A)
\end{equation}
and as a consequence of (\ref{idem-psiAV}) and (\ref{wmeas-wcp})  we have
\begin{equation}
\label{AV-axu2}
(\psi_{V}^{A}\otimes V)\circ (V\otimes (\psi_{V}^{A}\circ \psi_{A}^{V}))\circ
(\psi_{A}^{V}\otimes V)
\end{equation}
$$=(\mu_{A}\otimes V\otimes V)\circ (A\otimes \psi_{V}^{A}\otimes V)\circ (A\otimes
V\otimes (\psi_{V}^{A}\circ (V\otimes \eta_{A}))).$$

We now proceed to prove the twisted condition for $ \psi_{A}^{V}$ and
$\sigma_{A}^{V}$.

\begin{itemize}
\item[ ]$ \hspace{0.38cm} (V\otimes \mu_{A})\circ (\sigma_{A}^{V}\otimes A)\circ (V\otimes
\psi_{A}^{V})\circ (\psi_{A}^{V}\otimes V)$

\item[ ]$=\psi_{A}^{V}\circ (\mu_{A}\otimes V)\circ (A\otimes \psi_{V}^{A})\circ
(\sigma_{V}^{A}\otimes A)\circ (V\otimes \psi_{A}^{V})\circ (\psi_{A}^{V}\otimes V)
$

\item[ ]$=\psi_{A}^{V}\circ (\mu_{A}\otimes V)\circ (A\otimes \sigma_{V}^{A})\circ
(\psi_{V}^{A}\otimes V)\circ
(V\otimes (\psi_{V}^{A}\circ \psi_{A}^{V}))\circ  (\psi_{A}^{V}\otimes V)  $

\item[ ]$= \psi_{A}^{V}\circ (\mu_{A}\otimes V)\circ (\mu_{A}\otimes \sigma_{V}^{A})\circ (A\otimes \psi_{V}^{A}\otimes V)\circ (A\otimes
V\otimes (\psi_{V}^{A}\circ (V\otimes \eta_{A})))
  $

\item[ ]$=\psi_{A}^{V}\circ (\mu_{A}\otimes V)\circ (A\otimes ((\mu_{A}\otimes V)\circ (A\otimes \sigma_{V}^{A})\circ   (\psi_{V}^{A}\otimes V)\circ (
V\otimes (\psi_{V}^{A}\circ (V\otimes \eta_{A}))))) $

\item[ ]$= \psi_{A}^{V}\circ (\mu_{A}\otimes V)\circ (A\otimes ((\mu_{A}\otimes V)\circ (A\otimes \psi_{V}^{A})\circ (\sigma_{V}^{A}\otimes \eta_{A}))) $

\item[ ]$=\psi_{A}^{V}\circ (\mu_{A}\otimes V)\circ (A\otimes (\nabla_{A\otimes V}\circ \sigma_{V}^{A}))  $

\item[ ]$=\psi_{A}^{V}\circ (\mu_{A}\otimes V)\circ (A\otimes \sigma_{V}^{A}) $

\item[ ]$=(V\otimes \mu_{A})\circ (\psi_{A}^{V}\otimes A)\circ (A\otimes \sigma_{A}^{V}). $

\end{itemize}

The first equality follows by (\ref{AV-axu1}), the second and the fourth ones by the twisted condition for $(A\otimes V, \mu_{A\otimes V})$, the third one by (\ref{AV-axu2}), the fifth one by the associativity of $\mu_{A}$, the sixth one by the definition of $\nabla_{A\otimes V}$, the seventh one by (\ref{idemp-sigma-inv}) and finally the eighth  one by (\ref{idem-sim}).

The proof for the cocycle condition is the following:

\begin{itemize}

\item[ ]$ \hspace{0.38cm}(V\otimes \mu_{A})\circ (\sigma_{A}^{V}\otimes A)\circ  (V\otimes \sigma_{A}^{V})$

\item[ ]$= (V\otimes \mu_{A})\circ ((\nabla_{V\otimes A}\circ \psi_{A}^{V})\otimes A)\circ (\sigma_{V}^{A}\otimes A)\circ
    (V\otimes (\psi_{A}^{V}\circ \sigma_{V}^{A})) $

\item[ ]$= \psi_{A}^{V}\circ \psi_{V}^{A}\circ (V\otimes \mu_{A})\circ ((\psi_{A}^{V}\circ \sigma_{V}^{A})\otimes A)\circ (V\otimes (\psi_{A}^{V}\circ \sigma_{V}^{A})) $

\item[ ]$=\psi_{A}^{V}\circ (\mu_{A}\otimes V)\circ (A\otimes \psi_{V}^{A})\circ (\sigma_{V}^{A}\otimes A)\circ (V\otimes (\psi_{A}^{V}\circ \sigma_{V}^{A}))   $

\item[ ]$= \psi_{A}^{V}\circ (\mu_{A}\otimes V)\circ (A\otimes \sigma_{V}^{A})\circ (\psi_{V}^{A}\otimes V)\circ (V\otimes \sigma_{V}^{A}) $

\item[ ]$= \psi_{A}^{V}\circ (\mu_{A}\otimes V)\circ (A\otimes \sigma_{V}^{A})\circ (\sigma_{V}^{A}\otimes V) $

\item[ ]$= \psi_{A}^{V}\circ (\mu_{A}\otimes V)\circ (A\otimes (\nabla_{A\otimes V}\circ \sigma_{V}^{A}))\circ (\sigma_{V}^{A}\otimes V) $

\item[ ]$= \psi_{A}^{V}\circ (\mu_{A}\otimes V)\circ (A\otimes ((\mu_{A}\otimes V)\circ (A\otimes \psi_{V}^{A})\circ ( \sigma_{V}^{A}\otimes \eta_{A})))\circ (\sigma_{V}^{A}\otimes V) $

\item[ ]$=\psi_{A}^{V}\circ (\mu_{A}\otimes V)\circ (\mu_{A}\otimes \sigma_{V}^{A})\circ (A\otimes \psi_{V}^{A}\otimes V)\circ (\sigma_{V}^{A}\otimes (\psi_{V}^{A}\circ (V\otimes \eta_{A})) ) $

\item[ ]$= \psi_{A}^{V}\circ (\mu_{A}\otimes V)\circ (A\otimes \sigma_{V}^{A})\circ ( \psi_{V}^{A}\otimes V)\circ (V\otimes (\psi_{V}^{A}\circ \psi_{A}^{V}))\circ ((\psi_{A}^{V}\circ \sigma_{V}^{A})\otimes V)$

\item[ ]$=\psi_{A}^{V}\circ (\mu_{A}\otimes V)\circ (A\otimes \psi_{V}^{A})\circ ( \sigma_{V}^{A}\otimes V)\circ (V\otimes \psi_{A}^{V})\circ ((\psi_{A}^{V}\circ \sigma_{V}^{A})\otimes V)   $

\item[ ]$= (A\otimes \mu_{A})\circ ( \sigma_{A}^{V}\otimes A)\circ (V\otimes \psi_{A}^{V})\circ ( \sigma_{A}^{V}\otimes V).$

\end{itemize}

The first equality follows by $\nabla_{V\otimes A}\circ \psi_{A}^{V}=\psi_{A}^{V}$, the second one by (\ref{idem-psiVA}) and right $A$-linearity of $\nabla_{V\otimes A}$ where $\phi_{V\otimes A}=V\otimes \mu_{A}$. The third one is a consequence of (\ref{wmeas-wcp}) and (\ref{idemp-sigma-inv}), the forth one follows by the twisted condition for $(A\otimes V,\mu_{A\otimes V})$ as well as (\ref{idemp-sigma-inv}), the fifth one by the cocycle condition for $(A\otimes V,\mu_{A\otimes V})$, the sixth one by (\ref{idemp-sigma-inv}), the seventh one by the definition of $\nabla_{A\otimes V}$, the eighth one by the twisted condition for $(A\otimes V,\mu_{A\otimes V})$ and the associativity of $\mu_{A}$, the ninth one by (\ref{AV-axu2}), the tenth one by the twisted condition for $(A\otimes V,\mu_{A\otimes V})$ and the eleventh one by (\ref{AV-axu1}).

Therefore, $(V\otimes A,\mu_{V\otimes A})$ is a weak crossed product with associated
morphisms $\psi_{A}^{V}$, $\sigma_{A}^{V}$.

Finally, by a similar calculus we have $\mu_{V\otimes A} = \psi_{A}^{V}\circ \mu_{A\otimes V}\circ (\psi_{V}^{A}\otimes \psi_{V}^{A})$. Indeed:
\begin{itemize}

\item[ ]$ \hspace{0.38cm}\mu_{V\otimes A} $

\item[ ]$= (V\otimes \mu_{A})\circ (((V\otimes \mu_{A})\circ ((\psi_{A}^{V}\circ \sigma_{V}^{A})\otimes A)) \otimes A)\circ (V\otimes \psi_{A}^{V}\otimes A) $

\item[ ]$=\psi_{A}^{V}\circ (\mu_{A}\otimes V)\circ (A\otimes \psi_{V}^{A})\circ (\sigma_{V}^{A}\otimes \mu_{A})\circ  (V\otimes \psi_{A}^{V}\otimes A)  $

\item[ ]$= \psi_{A}^{V}\circ (\mu_{A}\otimes V)\circ (A\otimes \sigma_{V}^{A})\circ (\psi_{V}^{A}\otimes V)\circ (V\otimes (\psi_{V}^{A}\circ (V\otimes \mu_{A})\circ (\psi_{A}^{V}\otimes A)))  $

\item[ ]$= \psi_{A}^{V}\circ (\mu_{A}\otimes V)\circ (A\otimes \sigma_{V}^{A})\circ (\psi_{V}^{A}\otimes V)\circ (V\otimes ((\mu_{A}\otimes V)\circ (\mu_{A}\otimes \psi_{V}^{A})\circ $

\item[ ]$\hspace{0.38cm} (A\otimes (\psi_{V}^{A}\circ (V\otimes \eta_{A}))\otimes A)))$

\item[ ]$=\psi_{A}^{V}\circ (\mu_{A}\otimes V)\circ (\mu_{A}\otimes \sigma_{V}^{A})\circ (A\otimes \psi_{V}^{A}\otimes V)\circ (\psi_{V}^{A}\otimes \psi_{V}^{A})  $

\item[ ]$= \psi_{A}^{V}\circ \mu_{A\otimes V}\circ (\psi_{V}^{A}\otimes \psi_{V}^{A}). $
\end{itemize}

To complete the proof suppose that $(A\otimes V,\mu_{A\otimes V})$ is a weak crossed product with preunit $\nu:K\rightarrow A\otimes V$ and define $\upsilon=\psi_{A}^{V}\circ \nu:K\rightarrow V\otimes A$. To obtain that $\upsilon$ is a preunit for $(V\otimes A, \mu_{V\otimes A})$ it is sufficient to check:

\begin{equation}\label{sim-pre1-wcp}
    (V\otimes \mu_A)\circ (\sigma_{A}^{V}\otimes A)\circ
    (V\otimes \psi_{A}^{V})\circ (\upsilon\otimes V) =
    \nabla_{V\otimes A}\circ
    (V\otimes \eta_A),
    \end{equation}
    \begin{equation}\label{sim-pre2-wcp}
    (V\otimes \mu_A)\circ (\sigma_{A}^{V}\otimes A)\circ
    (V\otimes \upsilon)  =
    \nabla_{V\otimes A}\circ
    (V\otimes \eta_A),
    \end{equation}
    \begin{equation}\label{sim-pre3-wcp}
(V\otimes \mu_A)\circ (\psi_{A}^{V}\otimes A)\circ (A\otimes \upsilon)
= \beta_{\upsilon},
\end{equation}
i.e. the symmetric versions of the equalities (\ref{pre1-wcp}), (\ref{pre2-wcp}) and (\ref{pre3-wcp}).

Equality (\ref{sim-pre1-wcp}) follows by:

\begin{itemize}

\item[ ]$ \hspace{0.38cm} (V\otimes \mu_A)\circ (\sigma_{A}^{V}\otimes A)\circ
    (V\otimes \psi_{A}^{V})\circ (\upsilon\otimes V) $

\item[ ]$=(V\otimes\mu_{A})\circ ( \psi_{A}^{V}\otimes A)\circ (A\otimes (\psi_{A}^{V}\circ \sigma_{V}^{A}))\circ (\nu\otimes V)   $

\item[ ]$=\psi_{A}^{V}\circ (\mu_{A}\otimes V)\circ  (A\otimes \sigma_{V}^{A})\circ (\nu\otimes V)   $

\item[ ]$= \psi_{A}^{V}\circ \nabla_{A\otimes V}\circ (\eta_{A}\otimes V)  $

\item[ ]$= \psi_{A}^{V}\circ (\mu_{A}\otimes V)\circ (\eta_{A}\otimes (\psi_{V}^{A}\circ (V\otimes \eta_{A})))  $

\item[ ]$= \psi_{A}^{V}\circ \psi_{V}^{A}\circ (V\otimes \eta_{A}) $

\item[ ]$=\nabla_{V\otimes A}\circ (V\otimes \eta_{A}), $

\end{itemize}

where the first identity is consequence of the twisted condition for $(V\otimes A, \mu_{V\otimes A})$, the second one of (\ref{idem-sim}), the third one of (\ref{pre2-wcp}), the fourth one of the definition of $\nabla_{A\otimes V}$, the fifth one of the unit of $\mu_{A}$ and the sixth one of (\ref{idem-psiVA}).

The identity (\ref{sim-pre2-wcp}) follows by

\begin{itemize}

\item[ ]$ \hspace{0.38cm}\nabla_{V\otimes A}\circ
    (V\otimes \eta_A)  $

\item[ ]$=\psi_{A}^{V}\circ \nabla_{A\otimes V}\circ (\eta_{A}\otimes V)   $

\item[ ]$=\psi_{A}^{V}\circ  (\mu_{A}\otimes V)\circ (A\otimes \sigma_{V}^{A})\circ (\psi_{V}^{A}\otimes V)\circ (V\otimes \nu)   $

\item[ ]$= (V\otimes \mu_{A})\circ (\psi_{A}^{V}\otimes A)\circ (A\otimes ( \psi_{A}^{V}\circ \sigma_{V}^{A}))\circ (\psi_{V}^{A}\otimes V)\circ (V\otimes \nu)  $

\item[ ]$= (V\otimes \mu_{A})\circ  (( \psi_{A}^{V}\circ \sigma_{V}^{A})\otimes A)\circ (V\otimes \psi_{A}^{V})\circ (\nabla_{V\otimes A}\otimes V)\circ (V\otimes \nu)   $

\item[ ]$=(V\otimes \mu_{A})\circ  (( \psi_{A}^{V}\circ \sigma_{V}^{A})\otimes \mu_{A})\circ (V\otimes \psi_{A}^{V}\otimes A)\circ ((\psi_{A}^{V}\circ (\eta_{A}\otimes V))\otimes \upsilon)   $

\item[ ]$= (V\otimes \mu_{A})\circ  (( (V\otimes\mu_{A})\circ (\sigma_{A}^{V}\otimes A)\circ (V\otimes \psi_{A}^{V}))\otimes A)\circ ((\psi_{A}^{V}\circ (\eta_{A}\otimes V))\otimes \upsilon)  $

\item[ ]$=(V\otimes \mu_{A})\circ (\psi_{A}^{V}\otimes A)\circ (\mu_{A}\otimes V\otimes A)\circ (\eta_{A}\otimes ((\sigma_{V}^{A}\otimes A)\circ (V\otimes \upsilon)))   $

\item[ ]$=  (V\otimes \mu_A)\circ (\sigma_{A}^{V}\otimes A)\circ
    (V\otimes \upsilon), $

\end{itemize}

where the first equality follows by (\ref{idem-psiVA}), the second one by (\ref{pre1-wcp}), the third one by (\ref{idem-sim}), the fourth one by the twisted condition for $(V\otimes A, \mu_{V\otimes A})$ and (\ref{idem-psiVA}), the fifth one by the definition of $\nabla_{V\otimes A}$ and (\ref{idem-sim}), the sixth one  by the assocativity of $\mu_{A}$, the seventh one by the twisted condition for $(V\otimes A, \mu_{V\otimes A})$ and finally the eighth one by the properties of the unit of $A$.

Finally, by (\ref{pre3-wcp}) and similar arguments to the ones used previously, the proof for  (\ref{sim-pre3-wcp}) is the following:

\begin{itemize}

\item[ ]$ \hspace{0.38cm} (V\otimes \mu_{A})\circ (\psi_{A}^{V}\otimes A)\circ (A\otimes \upsilon)$

\item[ ]$=\psi_{A}^{V}\circ \beta_{\nu}   $

\item[ ]$= \psi_{A}^{V}\circ (\mu_{A}\otimes V)\circ (A\otimes \psi_{V}^{A})\circ (\nu\otimes A)  $

\item[ ]$=(V\otimes \mu_{A})\circ (\psi_{A}^{V}\otimes A)\circ (A\otimes \nabla_{V\otimes A})\circ (\nu\otimes A)    $

\item[ ]$=(V\otimes \mu_{A})\circ (((V\otimes \mu_{A})\circ (\psi_{A}^{V}\otimes A)\circ (A\otimes (\psi_{A}^{V}\circ (\eta_{A}\otimes V))))\otimes A)\circ  (\nu\otimes A)   $

\item[ ]$= (V\otimes \mu_{A})\circ (\upsilon \otimes A) $

\item[ ]$=\beta_{\upsilon}.  $

\end{itemize}

Therefore, $\upsilon$ is a preunit for $(V\otimes A, \mu_{V\otimes A})$ and the proof is finished.
\end{dem}

\begin{parrafo}
As our aim is to deal with algebraic structures that involve crossed products and coproducts, in the  rest of the section we give some definitions an results related to weak crossed coproducts. This theory is obtained by an straightforward dualization of the one of weak crossed products. The proofs can be obtained from the corresponding proofs for weak crossed
products defined in $A\otimes V$ or  $V\otimes A$ by dualization. We refer to
\cite{mra-preunit} for the notation and the main results.

Let $C$ be a coalgebra and $V$ an object. Suppose
that there exists a morphism $\chi_{V}^{C}:C\otimes V\rightarrow V\otimes C$
such that the following equality holds:
\begin{equation}\label{co-wmeas-wcp}
(\chi_{V}^{C}\otimes C)\circ (C\otimes \chi_{V}^{C})\circ
(\delta_{C}\otimes V)=(V\otimes \delta_{C})\circ \chi_{V}^{C}
\end{equation}

As a consequence the morphism $\Gamma_{C\otimes V}:C\otimes V\rightarrow C\otimes V$ defined by
\begin{equation}\label{co-idem-wcp}
\Gamma_{C\otimes V}=(C\otimes V\otimes \varepsilon_{C})\circ
(C\otimes \chi_{V}^{C})\circ (\delta_{C}\otimes V).
\end{equation}
 is  idempotent. Moreover, $\Gamma_{C\otimes V}$ is a  left $C$-comodule morphism
where the left coaction is given as $\rho_{C\otimes V}=\delta_{C}\otimes V$.

Henceforth we will consider quadruples $(C, V, \chi_{V}^{C}, \tau_{V}^{C})$ where $C$,
$V$ and $\chi_{V}^{C}$ satisfy the condition (\ref{co-wmeas-wcp}) and
$\tau_{V}^{C}:C\otimes V\rightarrow V\otimes V$ is a morphism in ${\mathcal C}$.
For the morphism $\Gamma_{C\otimes V}$ defined in (\ref{co-idem-wcp}) we denote by $C\Box V$ the image of $\Gamma_{C\otimes V}$ and by
$i_{C\otimes V}:C\Box V\rightarrow C\otimes V$, $p_{C\otimes
V}:C\otimes V\rightarrow C\Box V$ the injection and the projection
associated to the idempotent.

Following the ideas behind the theory of weak crossed products, we will set two properties that guarantee the coassociativity of certain weak crossed coproduct on $C\otimes V$.
We will say that a quadruple $(C, V, \chi_{V}^{C}, \tau_{V}^{C})$  satisfies the
twisted
condition if
\begin{equation}\label{co-twis-wcp}
(\tau_{V}^{C}\otimes C)\circ (C\otimes \chi_{V}^{C})\circ
(\delta_{C}\otimes V)=(V\otimes \chi_{V}^{C})\circ
(\chi_{V}^{C}\otimes V)\circ (C\otimes \tau_{V}^{C})\circ
(\delta_{C}\otimes V),
\end{equation}
 and the cocycle
condition holds if
\begin{equation}\label{co-cocy-wcp}
(\tau_{V}^{C}\otimes V)\circ (V\otimes \tau_{V}^{C})\circ
(\delta_{C}\otimes V)
\end{equation}
$$=(V\otimes \tau_{V}^{C})\circ (\chi_{V}^{C}\otimes V)\circ
(C\otimes \tau_{V}^{C})\circ (\delta_{C}\otimes V)
.$$

By virtue of Proposition 5.9  of \cite{mra-preunit} we will consider from now on,
and without loss of generality, that
\begin{equation}
\tau_{V}^{C}\circ\Gamma_{C\otimes V} = \tau_{V}^{C}
\end{equation}
for all quadruples $(C, V, \chi_{V}^{C}, \tau_{V}^{C})$.

For a quadruple $(C, V, \chi_{V}^{C}, \tau_{V}^{C})$  define the coproduct 
\begin{equation}\label{co-prod-todo-wcp}
\delta_{C\otimes  V} =(C\otimes \chi_{V}^{C}\otimes V)\circ
(\delta_{C}\otimes \tau_{V}^{C})\circ (\delta_{C}\otimes V)
\end{equation}
and let $\delta_{C\Box V}$ be the coproduct 
\begin{equation}
\label{co-prod-wcp} \delta_{C\Box V} = (p_{C\otimes V}\otimes
p_{C\otimes V})\circ\delta_{C\otimes V}\circ i_{C\otimes V}.
\end{equation}

If the twisted condition (\ref{co-twis-wcp}) and
the cocycle condition (\ref{co-cocy-wcp}) hold, $\delta_{C\otimes V}$ is a coassociative
coproduct that it is normalized with respect to $\Gamma_{C\otimes
V}$ (i.e. $\delta_{C\otimes V}\circ \Gamma_{C\otimes
V}=\delta_{C\otimes V}=(\Gamma_{C\otimes
V}\otimes \Gamma_{C\otimes
V})\circ \delta_{C\otimes V}$). As a consequence $\delta_{C\Box V}$ is also a
coassociative coproduct (Proposition 5.10 of \cite{mra-preunit}). Under these circumstances we say that $(C\otimes V, \delta_{C\otimes V})$ is a weak crossed coproduct.

Let $C$ be a coalgebra and $V$ and object in ${\mathcal C}$. If $\delta_{C\otimes
V}$ is a coassociative coproduct defined in
$C\otimes V$ with  precounit $\upsilon: C\otimes V\rightarrow K$ , i.e. a morphism
satisfying
\begin{equation}
(C\otimes V\otimes \upsilon)\circ \delta_{C\otimes V}=(\upsilon\otimes C\otimes
V)\circ \delta_{C\otimes V}=(C\otimes V\otimes ((\upsilon\otimes \upsilon)\circ
\delta_{C\otimes V}))\circ \delta_{C\otimes V},
\end{equation}
we obtain that $C\Box V$ is a coalgebra with coproduct
$$\delta_{C\Box V}=(p_{C\otimes V}^{\upsilon}\otimes p_{C\otimes V}^{\upsilon})\circ
\delta_{C\otimes V}\circ i_{C\otimes V}^{\upsilon}$$
and counit $\varepsilon_{C\Box V}=\upsilon\circ i_{C\otimes V}^{\upsilon}$, where
$p_{C\otimes V}^{\upsilon}$ and $i_{C\otimes V}^{\upsilon}$ are the injection and
the projection associated to the idempotent $\Gamma_{C\otimes
V}^{\upsilon}=(C\otimes V\otimes \upsilon)\circ \delta_{C\otimes V}:C\otimes
V\rightarrow C\otimes V$ (see Proposition 5.5 of \cite{mra-preunit}).

If moreover, $\delta_{C\otimes V}$ is left $C$-colinear for the coactions
$\rho_{C\otimes V}=\delta_{C}\otimes V$ and $\rho_{C\otimes V\otimes C\otimes V}=\rho_{C\otimes
V}\otimes C\otimes V$  and normalized with respect to $\Gamma_{C\otimes
V}^{\upsilon}$, the morphism
\begin{equation}
\label{gamma-nu}
\gamma_{\upsilon}:C\otimes V\rightarrow  C,\; \gamma_{\upsilon} =
(C\otimes \upsilon)\circ (\delta_{C}\otimes V)
\end{equation}
is comultiplicative and left $C$-colinear for $\rho_{C}=\delta_{C}$. Although
$\gamma_{\upsilon}$
is not a coalgebra morphism, because $C\otimes V$ is not a
coalgebra, we have that $\varepsilon_{C}\circ \gamma_{\upsilon} =
\upsilon$, and as a consequence the morphism $\bar{\gamma_{\upsilon}}
=\gamma_{\upsilon}\circ  i_{C\otimes V}^{\upsilon}
:C\Box V\rightarrow C$ is a coalgebra morphism.
\end{parrafo}

In the following theorem we give a characterization of weak
crossed coproducts with precounit:

\begin{teorema}
\label{co-thm1-wcp} Let $C$ be a coalgebra, $V$ an object and
$\Delta_{C\otimes V}:C\otimes V\rightarrow  C\otimes V\otimes C\otimes V$ a
morphism of left $C$-comodules.

Then the following statements are equivalent:
\begin{itemize}
\item[(i)] The coproduct $\Delta_{C\otimes V}$ is coassociative with precounit
$\upsilon$ and normalized with respect to $\Gamma_{C\otimes V}^{\upsilon}.$

\item[(ii)] There exist morphisms $\chi_{V}^{C}:C\otimes V\rightarrow V\otimes
C$, $\tau_{V}^{C}:C\otimes V\rightarrow V\otimes V$ and $\upsilon:C\otimes
V\rightarrow K$ such that if $\delta_{C\otimes V}$ is the coproduct defined
in (\ref{co-prod-todo-wcp}), the pair $(C\otimes V, \delta_{C\otimes
V})$ is a weak crossed coproduct with $\Delta_{C\otimes V} =
\delta_{C\otimes V}$ satisfying:
    \begin{equation}\label{co-pre1-wcp}
    (V\otimes \upsilon)\circ (\chi_{V}^{C}\otimes V)\circ (C\otimes
\tau_{V}^{C})\circ (\delta_{C}\otimes V)=(\varepsilon_{C}\otimes V)\circ
\Gamma_{C\otimes V}
    \end{equation}
    \begin{equation}\label{co-pre2-wcp}
     (\upsilon\otimes V)\circ (C\otimes \tau_{C}^{V})\circ (\delta_{C}\otimes
V)=(\varepsilon_{C}\otimes V)\circ \Gamma_{C\otimes V}
    \end{equation}
    \begin{equation}\label{co-pre3-wcp}
(\upsilon\otimes C)\circ (C\otimes \chi_{V}^{C})\circ
(\delta_{C}\otimes V) = \gamma_{\upsilon},
\end{equation}
\end{itemize}
where $\gamma_{\upsilon}$ is the morphism defined in
(\ref{gamma-nu}). In this case $\upsilon$ is a precounit for
$\delta_{C\otimes V}$, the idempotent morphism of the weak crossed
coproduct $\Gamma_{C\otimes V}$ is the idempotent $\Gamma_{C\otimes
V}^{\nu}$, and we say that $(C\otimes V, \delta_{C\otimes V})$ is a
 weak crossed coproduct with precounit $\upsilon$.
\end{teorema}

As a corollary of the previous result we have:

\begin{corol}\label{co-corol-wcp}
If $(C\otimes V, \delta_{C\otimes V})$ is a weak crossed coproduct
with precounit $\upsilon$, then $C\Box V$ is a coalgebra with the
coproduct defined in (\ref{co-prod-wcp}) and  counit
$\varepsilon_{C\Box V}=\upsilon\circ i_{C\otimes V}$.
\end{corol}

Like weak crossed products, weak crossed coproducts are universal constructions too.

\begin{teorema}
\label{co-uni-1}
Let $D$ be a coalgebra in ${\mathcal C}$. Then it holds equivalently:

\begin{itemize}

\item[(i)] There exists a weak crossed coproduct $(C\otimes V, \delta_{C\otimes V})$ with precounit $\upsilon$ and an isomorphism of coalgebras $\varpi:D\rightarrow C\Box V$.

\item[(ii)] There exist a coalgebra $C$, an object $V$, morphisms $$p_{C}:D\rightarrow
C,\;\;\;\;p_{V}:D\rightarrow V,\;\;\;\;\Gamma_{C\otimes V}: C\otimes V\rightarrow C\otimes
V,\;\;\;\;\varpi:D\rightarrow C\Box V$$ such that $p_{C}$ is a coalgebra morphism, $\Gamma_{C\otimes V}$ is an idempotent
morphism  of left $C$-comodules for the coaction $\rho_{C\otimes V}=\delta_{C}\otimes V$ and
$\varpi$ is an isomorphism such that
$$i_{C\otimes V}\circ \varpi= (p_{C}\otimes p_{V})\circ \delta_{D}$$  where $C\Box V$ is the image of
$\Gamma_{C\otimes V}$ and $i_{C\otimes V}$ is the associated injection.

\item[(iii)] There exist a coalgebra $C$, an object $V$ and morphisms $p_{C}:D\rightarrow C$, $p_{V}:D\rightarrow V$ and $\hat{\varpi}:C\otimes V\rightarrow D$ that satisfy
\begin{itemize}
\item[(iii-1)] $p_{C}$ is a coalgebra morphism.

\item[(iii-2)] $\hat{\varpi}$ is a morphism of left $C$-comodules for $\rho_{C\otimes V}=\delta_{C}\otimes V$ and $\rho_{D}=(p_{C}\otimes D)\circ \delta_{D}$.

\item[(iii-3)] $\hat{\varpi}\circ (p_{C}\otimes p_{V})\circ \delta_{D}=id_{D}.$

\end{itemize}

\end{itemize}

\end{teorema}

\begin{nota}
\label{co-proyeccion}
Note that, if the conditions  (ii) or (iii) of Theorem \ref{co-uni-1} hold it is possible
to characterize the isomorphism $\varpi$ as
$$\varpi=p_{C\otimes V}\circ (p_{C}\otimes p_{V})\circ \delta_{D}.$$

Also, as in the crossed product case, the explicit expression for the morphisms $p_{C}$ and $p_{V}$ are:

$$p_{C}= \gamma_{\upsilon}\circ i_{C\otimes V}\circ \varpi  ,\;\;\;  p_{V}=(\varepsilon_{C}\otimes V)\circ i_{C\otimes V}\circ \varpi $$

where $\varpi$ is the isomorphism of coalgebras.

\end{nota}

As we have done for weak crossed products, we can obtain a weak crossed coproduct as a universal construction, that is:

\begin{teorema}
\label{co-universal}
Let $(C\otimes V, \delta_{C\otimes V})$ be a weak crossed coproduct  with precounit
$\nu$, and  $D$ be a coalgebra. Suppose that there exists morphisms
$p_{C}:D\rightarrow C$ and $p_{V}:D\rightarrow V$ such that $p_{C}$ is a coalgebra. Then the following assertions are equivalent:

\begin{itemize}

\item[(i)] Then there exists an unique coalgebra morphism $\varpi:D\rightarrow C\Box V$
such that $\gamma_{\upsilon}\circ i_{C\otimes V}\circ \varpi=p_{C}$ and
$(\varepsilon_{C}\otimes V)\circ i_{C\otimes V}\circ \varpi=p_{V}$.

\item[(ii)] The equalities
\begin{itemize}
\item[(ii-1)] $\upsilon \circ  (p_{C}\otimes p_{V})\circ
\delta_{D}=\varepsilon_{D},$
\item[(ii-2)] $\chi_{V}^{C}\circ (p_{C}\otimes p_{V})\circ \delta_{D}=(p_{V}\otimes
p_{C})\circ \delta_{D},$

\item[(ii-3)] $\tau_{V}^{C}\circ (p_{C}\otimes p_{V})\circ \delta_{D}=(p_{V}\otimes
p_{V})\circ \delta_{D}$

\end{itemize}
hold.
\end{itemize}
\end{teorema}

\begin{nota}\label{VC}
Following the theory of weak crossed products, now we are interested in obtaining a weak crossed coproduct not on $C\otimes V$ but on $V\otimes C$. In this case we must consider morphisms
$\chi_{C}^{V}:V\otimes C\rightarrow C\otimes V$,
$\tau_{C}^{V}:V\otimes C\rightarrow V\otimes V$ and the associated
idempotent, $\Gamma_{V\otimes C}:V\otimes C\rightarrow V\otimes C$,
is defined by
\begin{equation}
\label{co-otra-nabla}
\Gamma_{V\otimes C}=(\varepsilon_{C}\otimes V\otimes C)\circ (\chi_{C}^{V}\otimes C)\circ (V\otimes \delta_{C}).
\end{equation}
The induced coproduct $\delta_{V\otimes C}:V\otimes C\rightarrow
V\otimes C\otimes V\otimes C$ is
\begin{equation}
\label{co-otro-prod}
\delta_{V\otimes C}=(V\otimes \chi_{C}^{V}\otimes C)\circ (\tau_{C}^{V}\otimes \delta_{C})\circ (V\otimes \delta_{C})
\end{equation}
and the precounit is a morphism $\nu:V\otimes C\rightarrow K$.
\end{nota}

The following result is the dual of Proposition \ref{pro-copro} and shows that if there exists a weak crossed coproduct on  $C\otimes
V$ then it is possible to find one on $V\otimes C$ under certain conditions.

\begin{prop}
\label{co-pro-copro}
Let $(C\otimes V,\delta_{C\otimes V})$ be a weak crossed coproduct with associated
morphisms $\chi_{V}^{C}:C\otimes V\rightarrow V\otimes C$, $\tau_{V}^{C}:C\otimes
V\rightarrow V\otimes V$. Let
$(C,V, \chi_{C}^{V}:V\otimes C\rightarrow C\otimes V)$ be a triple satisfying
\begin{equation}
\label{co-idem-sim}
(C\otimes \chi_{C}^{V})\circ (\chi_{C}^{V}\otimes C)\circ (V\otimes \delta_{C})=(\delta_{C}\otimes V)\circ \chi_{C}^{V}.
\end{equation}

Suppose that $\Gamma_{C\otimes V}$ is the idempotent associated to
 $(C\otimes V,\delta_{C\otimes V})$ and $\Gamma_{V\otimes C}$ is the idempotent given by (\ref{co-otra-nabla}). Then if
\begin{equation}
\label{co-idem-psiAV}
\chi_{C}^{V}\circ \chi_{V}^{C}=\Gamma_{C\otimes V}
\end{equation}
and
\begin{equation}
\label{co-idem-psiVA}
\chi_{V}^{C}\circ \chi_{C}^{V}=\Gamma_{V\otimes C},
\end{equation}
the pair $(V\otimes C,\delta_{V\otimes C})$ is a weak crossed coproduct with associated
morphisms $\chi_{C}^{V}$, $\tau_{C}^{V}=\sigma_{V}^{C}\circ \chi_{C}^{V}:V\otimes
C\rightarrow V\otimes V$ and the coproduct $\delta_{V\otimes C}$ defined in
(\ref{co-otro-prod}) satisfies
\begin{equation}
\label{co-idem-pro}
\delta_{V\otimes C}=(\chi_{V}^{C}\otimes \chi_{V}^{C})\circ \delta_{C\otimes V}\circ \chi_{C}^{V}.
\end{equation}

Moreover, if $(C\otimes V,\delta_{C\otimes V})$ is a weak crossed coproduct with precounit
$\upsilon:C\otimes V\rightarrow K$, the pair $(V\otimes C,\delta_{V\otimes C})$ is a weak
crossed coproduct with precounit $\nu=\upsilon\circ \chi_{C}^{V}:
V\otimes C\rightarrow K$.
\end{prop}

\section{Weak crossed biproducts}
The aim of this section is to obtain a weak bialgebra in a braided category $\mathcal C$ whose product is given as a weak crossed product and whose coproduct is given as a weak crossed coproduct. That is, we give the notion of a weak crossed biproduct, and furthermore we characterize such objects. Hence in this section the category ${\mathcal C}$ will be braided with braid $c$. The
notions of weak bialgebra and weak Hopf algebra in a braided category were introduced
in \cite{NikaRamon4} and the main properties of these algebraic objects were obtained
in \cite{NikaRamon5}. Weak Hopf algebras
in a braided category were recently called weak Hopf monoids by
Pastro and Street (see \cite{PS}). In this last reference, the authors showed that
it is possible to obtain examples of weak  Hopf
algebras and weak  bialgebras in a braided setting working with separable
Frobenius algebras in ${\mathcal C}$.

The definition of weak bialgebra in a braided category $\mathcal C$ is the following:

\begin{definicion}
\label{weakbialg}
A weak
bialgebra in a braided category ${\mathcal C}$ with braiding $c$,
is an object in ${\mathcal C}$ with an algebra structure $(D,
\eta_{D},\mu_{D})$ and a coalgebra structure $(D,
\varepsilon_{D},\delta_{D})$ satisfying:
\begin{itemize}
\item[(i)] $\delta_{D}\circ \mu_{D}=(\mu_{D}\otimes \mu_{D})\circ (D\otimes
c_{D,D}\otimes D)\circ (\delta_{D}\otimes \delta_{D}).$

\item[(ii)]$\varepsilon_{D}\circ \mu_{D}\circ (\mu_{D}\otimes
D)=((\varepsilon_{D}\circ \mu_{D})\otimes (\varepsilon_{D}\circ \mu_{D}))\circ
(D\otimes \delta_{D}\otimes D)$ \item[ ]$=((\varepsilon_{D}\circ \mu_{D})\otimes
(\varepsilon_{D}\circ \mu_{D}))\circ (D\otimes
(c_{D,D}^{-1}\circ\delta_{D})\otimes D).$

\item[(iii)]$(\delta_{D}\otimes D)\circ \delta_{D}\circ \eta_{D}=(D\otimes
\mu_{D}\otimes D)\circ ((\delta_{D}\circ \eta_{D}) \otimes (\delta_{D}\circ
\eta_{D}))$  \item[ ]$=(D\otimes (\mu_{D}\circ c_{D,D}^{-1})\otimes D)\circ
((\delta_{D}\circ \eta_{D}) \otimes (\delta_{D}\circ \eta_{D})).$

\end{itemize}

If moreover, the following conditions hold,

\begin{itemize}
\item[(iv)] There exists a morphism $\lambda_{D}:D\rightarrow D$
in ${\mathcal C}$ (called the antipode of $D$) satisfying:
\begin{itemize}
\item[(iv-1)] $id_{D}\wedge \lambda_{D}=((\varepsilon_{D}\circ
\mu_{D})\otimes D)\circ (D\otimes c_{D,D})\circ ((\delta_{D}\circ \eta_{D})\otimes
D),$ \item[(iv-2)] $\lambda_{D}\wedge
id_{D}=(D\otimes(\varepsilon_{D}\circ \mu_{D}))\circ (c_{D,D}\otimes D)\circ
(D\otimes (\delta_{D}\circ \eta_{D})),$ \item[(iv-3)]$\lambda_{D}\wedge
id_{D}\wedge \lambda_{D}=\lambda_{D},$
\end{itemize}
\end{itemize}
 the weak bialgebra $D$ is a weak Hopf algebra in the braided category $\mathcal C$.
\end{definicion}

If $D$ is a weak bialgebra it is possible to define the endomorphisms of $D$, $\Pi_{D}^{L}$
(target morphism), $\Pi_{D}^{R}$ (source morphism),
$\overline{\Pi}_{D}^{L}$ and $\overline{\Pi}_{D}^{R}$ by
$\Pi_{D}^{L}=((\varepsilon_{D}\circ \mu_{D})\otimes D)\circ
(D\otimes c_{D,D})\circ ((\delta_{D}\circ \eta_{D})\otimes D),\;$
$\Pi_{D}^{R}=(D\otimes(\varepsilon_{D}\circ \mu_{D}))\circ
(c_{D,D}\otimes D)\circ (D\otimes (\delta_{D}\circ \eta_{D})),\;$
$\overline{\Pi}_{D}^{L}=(D\otimes (\varepsilon_{D}\circ
\mu_{D}))\circ ((\delta_{D}\circ \eta_{D})\otimes D),\;$ and
$\overline{\Pi}_{D}^{R}=((\varepsilon_{D}\circ \mu_{D})\otimes
D)\circ(D\otimes (\delta_{D}\circ \eta_{D})).$ It is straightforward
to show that they are idempotent (Proposition 2.9 of \cite{NikaRamon5}), and by Proposition 2.10 of \cite{NikaRamon5} the following identities hold
\begin{equation}
\label{pi1}
\Pi_{H}^{L}\circ
\overline{\Pi}_{D}^{L}=\Pi_{D}^{L},\;\;\;\Pi_{D}^{L}\circ
\overline{\Pi}_{D}^{R}=\overline{\Pi}_{D}^{R},\;\;\;\overline{\Pi}_{D}^{L}\circ
\Pi_{D}^{L}=\overline{\Pi}_{D}^{L},\;\;\;\overline{\Pi}_{D}^{R}\circ
\Pi_{D}^{L}=\Pi_{D}^{L},
\end{equation}
$$
\;\;\;\;\;\;\Pi_{D}^{R}\circ
\overline{\Pi}_{D}^{L}=\overline{\Pi}_{D}^{L},\;\;\;
\Pi_{D}^{R}\circ
\overline{\Pi}_{D}^{R}=\Pi_{D}^{R},\;\;\;\overline{\Pi}_{D}^{L}\circ
\Pi_{D}^{R}=\Pi_{D}^{R},\;\;\; \overline{\Pi}_{D}^{R}\circ
\Pi_{D}^{R}=\overline{\Pi}_{D}^{R}.
$$
If $D$ is a weak Hopf algebra the antipode satisfies
\begin{equation}
\label{pi2}
\Pi_{D}^{L}\circ \lambda_{D}=\Pi_{D}^{L}\circ
\Pi_{D}^{R}= \lambda_{D}\circ \Pi_{D}^{R},\;\;\;\Pi_{D}^{R}\circ
\lambda_{D}=\Pi_{D}^{R}\circ \Pi_{D}^{L}= \lambda_{D}\circ
\Pi_{D}^{L},
\end{equation}
$$\Pi_{D}^{L}=\overline{\Pi}_{D}^{R}\circ
\lambda_{D}=\lambda_{D}
\circ\overline{\Pi}_{D}^{L},\;\;\;\Pi_{D}^{R}=
\overline{\Pi}_{D}^{L}\circ \lambda_{D}=\lambda_{D} \circ
\overline{\Pi}_{D}^{R}.$$

Moreover, it is
antimultiplicative, anticomultiplicative and leaves the unit and the
counit invariant,i.e.:
\begin{equation}
\label{ant-1}
\lambda_{D}\circ \mu_{D}=\mu_{D}\circ c_{D,D}\circ (\lambda_{D}\otimes
\lambda_{D}),
\end{equation}
\begin{equation}
\label{ant-2}
\delta_{D}\circ \lambda_{D}=(\lambda_{D}\otimes \lambda_{D})\circ c_{D,D}\circ
\delta_{D},
\end{equation}
\begin{equation}
\label{ant-3}
\lambda_{D}\circ \eta_{D}=\eta_{D},\;\; \varepsilon_{D}\circ\lambda_{D}
=\varepsilon_{D}
\end{equation}
(see Proposition 2.20 of \cite{NikaRamon5}).

In the following definition we introduce the notion of weak crossed biproduct that
generalizes to the weak setting the definition of cross product bialgebra due to
Bespalov and Drabant \cite{bes-drab2}.

\begin{definicion}
\label{biproduct}
A weak bialgebra $D$ in a braided monoidal category ${\mathcal C}$ is called a weak crossed biproduct if there exist an algebra $A$, a
coalgebra $C$ and morphisms
$$\psi_{A}^{C}:A\otimes C\rightarrow C\otimes A, \;\;\;\sigma_{A}^{C}:C\otimes
C\rightarrow C\otimes A, \;\;\; \nu:K\rightarrow C\otimes A, $$
$$\chi_{A}^{C}:C\otimes A\rightarrow A\otimes C, \;\;\;\tau_{A}^{C}:C\otimes
A\rightarrow A\otimes A, \;\;\; \upsilon:C\otimes A\rightarrow K $$
such that
\begin{itemize}

\item[(i)] The
pair $(C\otimes A, \mu_{C\otimes A})$ is a weak crossed product with preunit $\nu$.

\item[(ii)] The pair $(C\otimes A, \delta_{C\otimes A})$ is a weak
crossed coproduct with precounit $\upsilon$.

\item[(iii)] $\nabla_{C\otimes A}=\Gamma_{C\otimes A}$ where $\nabla_{C\otimes A}$ is
the idempotent associated to the weak crossed product and $\Gamma_{C\otimes A}$ the
idempotent associated to the weak crossed coproduct.

\item[(iv)] There exists an isomorphism of algebras and coalgebras $\alpha:C\times
A\rightarrow D$ where $C\times A$ denotes the image of $\nabla_{C\otimes A}$ or,
equivalently, the image of $\Gamma_{C\otimes A}$ .

\item[(v)]  The preunit and the precounit satisfy
$$(\varepsilon_{C}\otimes A)\circ \nu=\eta_{A},$$
$$\upsilon\circ (C\otimes \eta_{A})=\varepsilon_{C},$$
respectively.

\end{itemize}
\end{definicion}

If we combine the symmetric version of Theorem \ref{uni-1} and Theorem
\ref{co-uni-1} we obtain a characterization of weak crossed biproducts as universal
constructions:

\begin{teorema}
\label{Teo-biproduct}
Let $D$ be a weak bialgebra in ${\mathcal C}$. The following  are equivalent.

\begin{itemize}

\item[(i)] $D$ is  a weak crossed biproduct.

\item[(ii)] There exists  morphisms $\pi_{D}:D\rightarrow D$ and $\theta_{D}:D\rightarrow D$ such that

\begin{itemize}

\item[(ii-1)] $\pi_{D}\circ \eta_{D}=\eta_{D}.$

\item[(ii-2)] $\mu_{D}\circ (\pi_{D}\otimes \pi_{D})=\pi_{D}\circ \mu_{D}\circ (\pi_{D}\otimes \pi_{D}).$

\item[(ii-3)] $\varepsilon_{D}\circ \theta_{D}=\varepsilon_{D}.$

\item[(ii-4)] $(\theta_{D}\otimes \theta_{D})\circ \delta_{D}= (\theta_{D}\otimes \theta_{D})\circ \delta_{D}\circ \theta_{D}.$

\item[(ii-5)] $(\theta_{D}\otimes \pi_{D})\circ \delta_{D}\circ \mu_{D}\circ (D\otimes \pi_{D})=(D\otimes \mu_{D})\circ (((\theta_{D}\otimes \pi_{D})\circ \delta_{D})\otimes \pi_{D})$.

\item[(ii-6)] $(\theta_{D}\otimes D)\circ \delta_{D}\circ \mu_{D}\circ (\theta_{D}\otimes \pi_{D})=(\theta_{D}\otimes (\mu_{D}\circ (\theta_{D}\otimes \pi_{D})))\circ (\delta_{D}\otimes D). $

\item[(ii-7)] $\theta_{D}\wedge \pi_{D}=id_{D}.$

\end{itemize}

\item[(iii)] The following assertions hold.

\begin{itemize}

\item[(iii-1)] There exist an algebra $A$ and morphisms $i_{A}:A\rightarrow D$, $p_{A}:D\rightarrow A$ such that $i_{A}$ is an algebra morphism and $p_{A}\circ i_{A}=id_{A}$.

\item[(iii-2)] There exist an coalgebra $C$ and morphisms $i_{C}:C\rightarrow D$, $p_{C}:D\rightarrow C$ such that $p_{C}$ is a coalgebra morphism and $p_{C}\circ i_{C}=id_{C}$.

\item[(iii-3)] There exists an idempotent morphism $\nabla_{C\otimes A}:C\otimes A\rightarrow C\otimes A$ of right $A$-modules, for the action $\phi_{C\otimes A}=C\otimes \mu_{A}$, and left $C$-comodules, for the coaction $\rho_{C\otimes A}=\delta_{C}\otimes A$, and an isomorphism $\omega:C\times A\rightarrow D$ such that
    $$p_{C\otimes A}=\omega^{-1}\circ \mu_{D}\circ (i_{C}\otimes i_{A})$$
    and
    $$i_{C\otimes A}=(p_{C}\otimes p_{A})\circ \delta_{D}\circ \omega$$
     where $C\times A$ denotes the image of $\nabla_{C\otimes A}$ and $p_{C\otimes A}$, $i_{C\otimes A}$ the associated projection and injection respectively.
\end{itemize}

\end{itemize}
\end{teorema}

\begin{dem}
(i) $\Rightarrow$ (iii) By the symmetric version of Theorem \ref{uni-1} we have that there exist morphisms $i_{A}:A\rightarrow
D$, $i_{C}:C\rightarrow D$, $\nabla_{C\otimes A}: C\otimes A\rightarrow A\otimes
C$ where $i_{A}$ is an algebra morphism, $\nabla_{C\otimes A}$ is an idempotent
morphism  of right $A$-modules for the action $\phi_{C\otimes A}=C\otimes \mu_{A}$ and $\alpha$ is an isomorphism of algebras such that
$\alpha\circ p_{C\otimes A}=\mu_{D}\circ (i_{C}\otimes i_{A})$. By (ii) of Definition \ref{biproduct} and Theorem \ref{co-uni-1} we obtain that $\nabla_{C\otimes A}$ is a morphism  of left $C$-comodules for the coaction $\rho_{C}=\delta_{C}\otimes V$ and there exists morphisms $p_{C}:D\rightarrow
C$, $p_{A}:D\rightarrow A$, where $p_{C}$ is a coalgebra morphism and
$i_{C\otimes A}\circ \alpha^{-1}=(p_{C}\otimes p_{A})\circ \delta_{D}$.
The morphisms $i_{A}$, $i_{C}$, $p_{C}$ and $p_{A}$ are defined by

\begin{equation}
\label{ipip}
i_{A}=\alpha\circ p_{C\otimes A}\circ \beta_{\nu},\;\;\;\; p_{A}=(\varepsilon_{C}\otimes A)\circ i_{C\otimes A}\circ \alpha^{-1},
\end{equation}
\begin{equation}
\label{ipip2}
i_{C}=\alpha\circ p_{C\otimes A}\circ (C\otimes \eta_{A}),\;\;\;\; p_{C}=\gamma_{\upsilon}\circ  i_{C\otimes A}\circ \alpha^{-1}.
\end{equation}

Then, using the condition of right $A$-module morphism for $\nabla_{C\otimes A}$ and (5) of Definition \ref{biproduct} we have:
$$p_{A}\circ i_{A}=(\varepsilon_{C}\otimes A)\circ \nabla_{C\otimes A}\circ (C\otimes \mu_{A})\circ (\nu\otimes A)=(\varepsilon_{C}\otimes \mu_{A})\circ ((\nabla_{C\otimes A}\circ \nu)\otimes A)=$$
$$(\varepsilon_{C}\otimes \mu_{A})\circ ( \nu\otimes A)=\mu_{A}\circ (\eta_{A}\otimes A)=id_{A}.$$

In a similar way, by the condition of left $C$-comodule morphism for $\nabla_{C\otimes A}$ and (5) of Definition \ref{biproduct}, we prove that $p_{C}\circ i_{C}=id_{C}.$

(iii) $\Rightarrow$ (i)  The proof of this implication is a direct consequence of Theorem \ref{co-uni-1} and
the symmetric version of Theorem \ref{uni-1}. Obviously, the isomorphism $\alpha:C\times A\rightarrow D$ is $\omega$ and the preunit and the precounit are:
\begin{equation}
\label{preunit-precounit}
\nu=i_{C\otimes A}\circ \omega^{-1}\circ \eta_{D},\;\;\; \upsilon=\varepsilon_{D}\circ \omega\circ p_{C\otimes A}.
\end{equation}
Then,  we have
$$(\varepsilon_{C}\otimes A)\circ \nu=(\varepsilon_{C}\otimes A)\circ i_{C\otimes A}\circ \omega^{-1}\circ \eta_{D}=(\varepsilon_{C}\otimes A)\circ (p_{C}\otimes p_{A})\circ \delta_{D}\circ \omega \circ \omega^{-1}\circ \eta_{D}=p_{A}\circ \eta_{D}=$$
$$p_{A}\circ i_{A}\circ \eta_{A}=\eta_{A},$$
and in a similar way we obtain $\upsilon\circ (C\otimes \eta_{A})=\varepsilon_{C}.$

(ii) $\Rightarrow$ (iii) First note that (ii-1) and (ii-2) imply that $\pi_{D}$ is an idempotent morphism. In a similar way $\theta_{D}$ is idempotent by  (ii-3) and (ii-4). Moreover, by (ii-7) it is easy to prove that
\begin{equation}\label{iden-d}
\nabla_{D\otimes D}=(\theta_{D}\otimes \pi_{D})\circ \delta_{D}\circ \mu_{D}\circ (\theta_{D}\otimes \pi_{D})
\end{equation}
is an idempotent morphism with splitting morphisms $$p_{D\otimes D}=\mu_{D}\circ (\theta_{D}\otimes \pi_{D}),\;\;\;i_{D\otimes D}=(\theta_{D}\otimes \pi_{D})\circ \delta_{D}.$$

We define $i_{A}$, $p_{A}$, $i_{C}$, $p_{C}$ as the morphisms  that split $\pi_{D}$ and $\theta_{D}$ respectively, i.e.,
$$i_{A}\circ p_{A}= \pi_{D},\;\;\;p_{A}\circ i_{A}= id_{A},$$
$$i_{C}\circ p_{C}= \theta_{D},\;\;\;p_{C}\circ i_{C}= id_{C},$$
for some objects $A$ and $C$. By (ii-1) and (ii-2) we obtain that $A$ is an algebra with unit $\eta_{A}=p_{A}\circ \eta_{D}$ and product
$\mu_{A}=p_{A}\circ  \mu_{D}\circ (i_{A}\otimes i_{A}).$
Similarly, by (ii-3) and (ii-4) we prove that $C$ is an coalgebra with counit $\varepsilon_{C}= \varepsilon_{D}\circ i_{C}$ and coproduct
$\delta_{C}=(p_{C}\otimes p_{C})\circ  \delta_{C}\circ i_{C}.$ Therefore, $i_{A}$ is an algebra morphism and $p_{C}$ a coalgebra morphism.

Now define $\nabla_{C\otimes A}:C\otimes A\rightarrow C\otimes A$ by
$$\nabla_{C\otimes A}=(p_{C}\otimes p_{A})\circ \nabla_{D\otimes D}\circ (i_{C}\otimes i_{A}),$$
where $\nabla_{D\otimes D}$ is the idempotent defined in (\ref{iden-d}).
Then, $\nabla_{C\otimes A}$ is idempotent, satisfies the identity
$$\nabla_{C\otimes A}=(p_{C}\otimes p_{A})\circ \delta_{D}\circ  \mu_{D}\circ (i_{C}\otimes i_{A})  $$
and, as a consequence,  the splitting morphisms are
$$p_{C\otimes A}=\mu_{D}\circ (i_{C}\otimes i_{A}),\;\;\; i_{C\otimes A}=(p_{C}\otimes p_{A})\circ \delta_{D},$$
the image of $\nabla_{C\otimes A}$ is $A\times C=D$ and $\omega=id_{D}=\varpi$.

Finally, $\nabla_{C\otimes A}$ is a morphism of right $A$-modules, for the action {\blue $\phi$} {\red $\varphi$}$_{C\otimes A}=C\otimes \mu_{A}$, and a morphism of left $C$-comodules, for the coaction $\rho_{C\otimes A}=\delta_{C}\otimes A$. Indeed, to obtain the $A$-linearity of $\nabla_{C\otimes A}$ check:

\begin{itemize}

\item[ ]$ \hspace{0.38cm} \nabla_{C\otimes A}\circ \phi_{C\otimes A} $

\item[ ]$=(p_{C}\otimes (p_{A}\circ \mu_{D}))\circ (D\otimes \pi_{D}\otimes D)\circ ((\delta_{D}\circ i_{C})\otimes (i_{A}\circ \mu_{A}))$

\item[ ]$= (p_{C}\otimes (p_{A}\circ \mu_{D}))\circ (D\otimes (\mu_{D}\circ (\pi_{D}\otimes D))\otimes D)\circ ((\delta_{D}\circ i_{C})\otimes i_{A}\otimes  i_{A}) $

\item[ ]$=(p_{C}\otimes (\mu_{A}\circ (p_{A}\otimes A)))\circ ((\nabla_{D\otimes D}\circ (i_{C}\otimes i_{A}))\otimes A) $

\item[ ]$=(C\otimes \mu_{A})\circ (\nabla_{C\otimes A}\otimes A),$

\end{itemize}
where the first equality follows by (ii-5), the second one by the associativity of $\mu_{D}$ and the third one by (ii-5).

Using (ii-6) and by similar computations, follows the proof for $\nabla_{C\otimes A}$ to be of left $C$-comodules.

(iii) $\Rightarrow$ (ii) Define the morphism $\pi_{D}:D\rightarrow D$ and $\theta_{D}:D\rightarrow D$ as $\pi_{D}=i_{A}\circ p_{A}$ and $\theta_{D}=i_{C}\circ p_{C}$. These morphisms are idempotent because  $p_{A}\circ i_{A}=id_{A}$ and  $p_{C}\circ i_{C}=id_{C}$.
To check (ii-1) compute:
$$\pi_{D}\circ \eta_{D}=i_{A}\circ p_{A}\circ \eta_{D}=i_{A}\circ \eta_{A}=\eta_{D}$$
where the second and the third identities follows by the condition of algebra morphism for $i_{A}$. Using as well that $i_A$ is an algebra morphism we obtain
$$\pi_{D}\circ \mu_{D}\circ (\pi_{D}\otimes \pi_{D})=i_{A}\circ \mu_{A}\circ (p_{A}\otimes p_{A})=\mu_{D}\circ (\pi_{D}\otimes \pi_{D}),$$
i.e., (ii-2) holds. In a similar way, using the condition of coalgebra morphism for $p_{C}$ we obtain (ii-3) and (ii-4).

Note that we have the following identities for the splitting morphisms associated to $\nabla_{C\otimes A}$,
$$p_{C\otimes A}=\varpi \circ \mu_{D}\circ (i_{C}\otimes i_{A}),\;\;\; i_{C\otimes A}=(p_{C}\otimes p_{A})\circ \delta_{D}\circ \omega.$$

Using these equalities compute:
\begin{itemize}

\item[ ]$ \hspace{0.38cm} \varpi\circ (\theta_{D}\wedge \pi_{D}) \circ \omega  $

\item[ ]$ =\varpi\circ \mu_{D}\circ (i_{C}\otimes i_{A})\circ (p_{C}\otimes p_{A})\circ \delta_{D}\circ \omega  $

\item[ ]$= p_{C\otimes A}\circ i_{C\otimes A}$

\item[ ]$=id_{C\times A}$

\end{itemize}
and therefore, $\theta_{D}\wedge \pi_{D}=id_{D}$, thus (ii-7) holds.

To prove (ii-5) we need some preliminary steps. First note that for the idempotent morphism $\nabla_{D\otimes D}= (\theta_{D}\otimes \pi_{D})\circ \delta_{D}\circ \mu_{D}\circ (\theta_{D}\otimes \pi_{D})$  the equality
$\nabla_{C\otimes A}=(p_{C}\otimes p_{A})\circ \nabla_{D\otimes D}\circ (i_{C}\otimes i_{A})$ holds and, using the condition of right $A$-module morphism for $\nabla_{C\otimes A}$:
\begin{itemize}

\item[ ]$ \hspace{0.38cm}  \nabla_{D\otimes D}\circ (D\otimes (\mu_{D}\circ (i_{A}\otimes i_{A}))) $

\item[ ]$ =(i_{C}\otimes i_{A})\circ \nabla_{C\otimes A}\circ (p_{C}\otimes \mu_{A}) $

\item[ ]$=(i_{C}\otimes (i_{A}\circ \mu_{A}))\circ ((\nabla_{C\otimes A}\circ (p_{C}\otimes A))\otimes A) $

\item[ ]$=(i_{C}\otimes (\mu_{D}\circ (i_{A}\otimes i_{A})))\circ  ((\nabla_{C\otimes A}\circ (p_{C}\otimes A))\otimes A) $

\item[ ]$ =(D\otimes \mu_{D})\circ ((\nabla_{D\otimes D}\circ (D\otimes i_{A}))\otimes i_{A})$

\end{itemize}
and then we obtain
\begin{equation}
\label{nabla-prin}
\nabla_{D\otimes D}\circ (D\otimes (\mu_{D}\circ (\pi_{D}\otimes \pi_{D})))=
(D\otimes \mu_{D})\circ ((\nabla_{D\otimes D}\circ (D\otimes \pi_{D}))\otimes \pi_{D}).
\end{equation}
As a consequence,

\begin{itemize}

\item[ ]$ \hspace{0.38cm}(\theta_{D}\otimes \pi_{D})\circ \delta_{D}\circ \mu_{D}\circ (D\otimes \pi_{D})$

\item[ ]$=(\theta_{D}\otimes \pi_{D})\circ \delta_{D}\circ \mu_{D}\circ ((\theta_{D}\wedge \pi_{D})\otimes \pi_{D})   $

\item[ ]$=\nabla_{D\otimes D}\circ (\theta_{D}\otimes (\mu_{D}\circ (\pi_{D}\otimes \pi_{D})))\circ (\delta_{D}\otimes D)  $

\item[ ]$=(D\otimes \mu_{D})\circ ((\nabla_{D\otimes D}\circ (\theta_{D}\otimes \pi_{D})\circ \delta_{D})\otimes \pi_{D}) $

\item[ ]$=(D\otimes \mu_{D})\circ (((\theta_{D}\otimes \pi_{D})\circ \delta_{D}\circ (\theta_{D}\wedge \pi_{D}))\otimes \pi_{D})  $

\item[ ]$=(D\otimes \mu_{D})\circ (((\theta_{D}\otimes \pi_{D})\circ \delta_{D})\otimes \pi_{D}),$

\end{itemize}
where the first equality follows by (ii-7),  the second one by the associativity of $\mu_{D}$ and (ii-2), the third one by (\ref{nabla-prin}), the fourth one by definition of $\nabla_{D\otimes D}$ and finally the last one by (ii-7) and these computations yield (ii-5).

The proof of (ii-6) is similar and we leave it to the reader.

\end{dem}

If we are under the conditions of Theorem \ref{Teo-biproduct}, observe that
equalities (ii-1) and (ii-2) mean that $i_A$ is a morphism of algebras. Similarly,
(ii-3) and (ii-4) mean that $p_C$ is a morphism of coalgebras. Now recall from the
symmetric version of Theorem \ref{uni-1} that we require the existence of a morphism
$\hat{\omega}:D\rightarrow C\otimes A$ that must be of right $A$-modules and must
satisfy $\mu_D\circ (i_C\otimes i_A)\circ \hat{\omega} = id_D$, and moreover
$\hat{\omega}\circ \mu_D\circ (i_C\otimes i_A) = \nabla_{C\otimes A}$. In this case
\[\hat{\omega} = (p_C\otimes p_A)\circ \delta_D\]
and condition (ii-5) guarantees that it is a morphism of right $A$-modules. In an
analogous way, and by Theorem \ref{co-uni-1}, there exists a morphism
$\hat{\varpi}:C\otimes A\rightarrow D$ that satisfies that $\hat{\varpi}\circ
(p_C\otimes p_A)\circ \delta_D = id_D$ and $(p_C\otimes p_A)\circ \delta_D\circ
\hat{\varpi} = \nabla_{C\otimes A}$. In this case, $\hat{\varpi}$ must be of left
$C$-comodules. This property is obtained by (ii-6) and moreover we have
\[\hat{\varpi} = \mu_D\circ (i_C\otimes i_V).\]
Thus, the meaning of (ii-5) and (ii-6) is that $\nabla_{C\otimes A}$ is a morphism
of right $A$-modules and of left $C$-comodules. Finally, equality (ii-7) assures
that $\nabla_{C\otimes A}$ is an idempotent morphism.

\begin{nota}
In the conditions of the previous result, if $\nabla_{C\otimes A}=id_{C\otimes A}$ we obtain that $\nabla_{D\otimes D}=\theta_{D}\otimes \pi_{D}$. Then, the universal characterization of cross product bialgebras  obtained by Bespalov and Drabant in Theorem 3.2 of  \cite{bes-drab2} is a particular instance of the universal property obtained for weak crossed biproducts in Theorem \ref{Teo-biproduct}.

\end{nota}

\section{Weak projections and weak crossed biproducts}

This section is devoted to study weak bialgebras with a weak projection living in a braided monoidal category $\mathcal C$. To do this we use the theory developed in Section 3, and we obtain that every weak bialgebra with a weak projection is a weak crossed biproduct whose comultiplication is a classical cosmash coproduct when it is particularized to the non weak case.

\begin{definicion}
\label{modulo-coalgebra}
Let $B$ be a weak Hopf algebra in ${\mathcal C}$ and let $(D, \phi_{D})$ be a
coalgebra which is also a right $B$-module and such that
\begin{equation}
\label{comp}
(\phi_{D}\otimes \phi_{D})\circ (D\otimes c_{D,B}\otimes B)\circ (\delta_{D}\otimes
\delta_{B})=\delta_{D}\circ \phi_{D}.
\end{equation}
We say that $D$ is a right $B$-module
coalgebra if the following equivalent conditions hold:
\begin{equation}
\label{i)}
\varepsilon_{D}\circ \phi_{D}\circ (D\otimes \mu_{B})=
(\varepsilon_{D}\otimes \varepsilon_{B})\circ (\phi_{D}\otimes
\mu_{B})\circ (D\otimes (c_{B,B}^{-1}\circ \delta_{B})\otimes B),
\end{equation}
\begin{equation}
\label{ii)}
\varepsilon_{D}\circ \phi_{D}\circ (D\otimes \mu_{B})=
(\varepsilon_{D}\otimes \varepsilon_{B})\circ (\phi_{D}\otimes
\mu_{B})\circ (D\otimes \delta_{B}\otimes B),
\end{equation}
\begin{equation}
\label{iii)}
\phi_{D}\circ (D\otimes \Pi_{B}^{L})=(D\otimes
(\varepsilon_{D}\circ \phi_{D}))\circ ((c_{D,D}^{-1}\circ
\delta_{D})\otimes B),
\end{equation}
\begin{equation}
\label{iv)}
\phi_{D}\circ (D\otimes
\overline{\Pi}_{B}^{L})=(D\otimes (\varepsilon_{D}\circ
\phi_{D}))\circ (\delta_{D}\otimes B),
\end{equation}
\begin{equation}
\label{v)}
\varepsilon_{D}\circ \phi_{D}\circ (D\otimes
\Pi_{B}^{L})=\varepsilon_{D}\circ \phi_{D},
\end{equation}
\begin{equation}
\label{vi)}
\varepsilon_{D}\circ \phi_{D}\circ (D\otimes
\overline{\Pi}_{B}^{L})=\varepsilon_{D}\circ \phi_{D}.
\end{equation}

Under these conditions it is  easy to prove that $(B,D, \psi_{RR}=(B\otimes
\phi_{D})\circ (c_{D,B}\otimes B)\circ (D\otimes \delta_{B}))$ is
a right-right weak entwining structure and $(D,\phi_{D},\delta_{D})\in {\mathcal
M}_{B}^{D}(\psi_{RR}).$

\end{definicion}

\begin{definicion}
\label{def-weakprojection}
Let $D$ be a weak bialgebra and $B$ a weak Hopf algebra in the category ${\mathcal C}$. Suppose that there exists a morphism $f:B\rightarrow D$ of weak bialgebras and a morphism $g:D\rightarrow B$ of coalgebras such that $g\circ f = id_B$. Define $\phi_D = \mu_D\circ (D\otimes f)$, that induces a right $B$-module structure on $D$, and let $g$ be a morphism of right $B$-modules, that is,
\begin{equation}
\label{bmodule}
g\circ \phi_{D}=\mu_{B}\circ (g\otimes B).
\end{equation}
 Under these circumstances we say that $D$ is a weak bialgebra with a weak projection onto $B$. Henceforth we will always use the notation $(D,B,f,g)$  to refer to a weak bialgebra with a weak projection.  In these conditions $g\circ \eta_{D}=\eta_{B}$.

If we consider the $B$-module structure on $D$ given by $\phi_D $, we obtain that $D$ is a right $B$-module coalgebra. Then $(B,D, \psi_{RR})$ with the morphism $\psi_{RR}:D\otimes B\rightarrow B\otimes D$ given by:
\begin{equation}
\label{psi}\psi_{RR} = (B\otimes \mu_D)\circ (c_{D,B}\otimes f)\circ (D\otimes \delta_B)
\end{equation}
is a right-right weak entwining structure such that $(D,\phi_{D}, \delta_{D})\in {\mathcal M}_B^D(\psi_{RR})$. Notice that the morphism $e_{RR}$ related to this weak entwining structure is given by:
\begin{equation}
\label{err-pi}
e_{RR} = \Pi_B^R\circ g.
\end{equation}

\end{definicion}

\begin{exemplo}
\label{example1}
As group algebras, that are the natural examples of
Hopf algebras, groupoid algebras  provide examples
of weak Hopf algebras. Recall that a groupoid $G$ is simply a
category in which every morphism is an isomorphism. In this
example, we consider finite groupoids, i.e. groupoids with a
finite number of objects. The set of objects of $G$, called also the base of $G$, will be
denoted by $G_{0}$ and the set of morphisms by $G_{1}$. The
identity morphism on $x\in G_{0}$ will  be denoted by $id_{x}$
and for a morphism $\sigma:x\rightarrow y$ in $G_{1}$, we write
$s(\sigma)$ and $t(\sigma)$, respectively for the source and the
target of $\sigma$.

Let $G$ be a groupoid and $R$ a commutative ring. The groupoid
algebra is the direct product in $R$-Mod
$$RG=\bigoplus_{\sigma\in G_{1}}R\sigma$$
with the product of two morphisms being equal to their composition
if the latter is defined and $0$ in otherwise, i.e.
$\mu_{RG}(\tau\otimes \sigma)=\tau\circ \sigma$ if $s(\tau)=t(\sigma)$ and
$\mu_{RG}(\tau\otimes \sigma)=0$ if $s(\tau)\neq t(\sigma)$. The unit element is
$1_{RG}=\sum_{x\in G_{0}}id_{x}$. The algebra $RG$ is a
cocommutative weak Hopf algebra, with coproduct $\delta_{RG}$,
counit $\varepsilon_{RG}$ and antipode $\lambda_{RG}$ given by the
formulas:
$$\delta_{RG}(\sigma)=\sigma\otimes \sigma, \;\;\;\varepsilon_{RG}
(\sigma)=1,\;\;\; \lambda_{RG}(\sigma)=\sigma^{-\dot{}1}.$$

For the weak Hopf algebra $RG$ the morphisms target and source are
respectively,
$$\Pi_{RG}^{L}(\sigma)=id_{t(\sigma)},\;\;\;
\Pi_{RG}^{R}(\sigma)=id_{s(\sigma)}$$ and $\lambda_{RG}\circ
\lambda_{RG}=id_{RG}$, i.e. the antipode is involutory.

A wide subgroupoid  of a groupoid $G$ is a groupoid $H$ provided with a functor $F:H\rightarrow G$ which is the identity on the objects, and induces inclusions $hom_{H}(x,y)\subset hom_{G}(x,y)$, i.e., it has the same base, and (perhaps) less arrows.

Let $G$ be a groupoid. An exact factorization of $G$ is a pair of wide subgroupoids of $G$, $H$ and $V$, such that for any $\sigma\in G_{1}$, there exist unique $\sigma_{V}\in V_{1}$, $\sigma_{H}\in H_{1}$, such that $\sigma=\sigma_{H}\circ \sigma_{V}$. Following the notation of \cite{MacK} and \cite{AN1} we denote $G$ as $H\bowtie V$ because in Theorems 2.10 and 2.15 of \cite{MacK} was proved  that the notion of groupoids with an exact factorization is equivalent to the notion of mathched pairs of groupoids  and to the notion of vacant double groupoid.

In this example, we will prove that any groupoid $G$ with an exact factorization $H\bowtie V$ induces a non-trivial example of a weak bialgebra with a weak projection. Put $B=RV$ and $D=RG$ and define
$f:B\rightarrow D$ by $f(\sigma)=\sigma$ and $g:D\rightarrow B$ by $g(\tau)=\tau_{V}$. Then, it is easy to show that $f$ is an algebra-coalgebra morphism and $g\circ f=id_{B}$. Moreover, $g$ is a coalgebra morphism because $\varepsilon_{B}\circ g=\varepsilon_{D}$ and
$$(\delta_{B}\circ g)(\tau)=\tau_{V}\otimes \tau_{V}=g(\tau)\otimes g(\tau)=((g\otimes g)\circ \delta_{D})(\tau).$$

The morphism $g$ is a right $B$-module morphism, i.e. satisfies (\ref{bmodule}), because
$$(g\circ \mu_{D}\circ (D\otimes f))(\tau\otimes \sigma)=\left\{ \begin{array}{ll}
0\;\;\;\;\;& {\rm if}\;  s(\tau)\neq t(\sigma) \\
\;\\
(\tau\circ \sigma)_{V} & {\rm if}\;  s(\tau)=t(\sigma) \end{array} \right. $$
and, on the other hand,
$$( \mu_{B}\circ (g\otimes B))(\tau\otimes \sigma)=\left\{ \begin{array}{ll}
0\;\;\;\;\;& {\rm if}\;  s(\tau_{V})\neq t(\sigma) \\
\;\\
\tau_{V}\circ \sigma & {\rm if}\;  s(\tau_{V})=t(\sigma) \end{array} \right. $$
But $s(\tau_{V})=s(\tau)$ and $(\tau\circ \sigma)_{V}=\tau_{V}\circ \sigma$ and then we have (\ref{bmodule}).

Finally, $g$ does not satisfy the conditon of algebra morphism. Indeed: let $\sigma$ be a morphism in $V_{1}$ different of the identities with $s(\sigma)\neq t(\sigma)$ and let $\tau$ a morphism in $H_{1} $ in the same conditions of $\sigma$ and such that $s(\tau)=t(\sigma)$. Put $\omega=\tau\circ \sigma$ and let $\gamma$ a morphism in $H_{1}$ such that $s(\gamma)=t(\tau)$. For example
\begin{center}
\scalebox{0.85}
{
\begin{pspicture}(0,-1.478125)(8.682813,1.478125)
\psline[linewidth=0.04cm,arrowsize=0.05291667cm 2.0,arrowlength=1.4,arrowinset=0.4]{->}(5.4409375,-0.9003125)(8.120937,-0.9403125)
\psline[linewidth=0.04cm,arrowsize=0.05291667cm 2.0,arrowlength=1.4,arrowinset=0.4]{->}(3.9409375,1.0996875)(5.1009374,-0.6603125)
\psline[linewidth=0.04cm,arrowsize=0.05291667cm 2.0,arrowlength=1.4,arrowinset=0.4]{->}(1.8409375,-0.5603125)(3.5209374,1.0596875)
\usefont{T1}{ptm}{m}{n}
\rput(1.5823437,-0.9103125){$x$}
\usefont{T1}{ptm}{m}{n}
\rput(5.052344,-0.9503125){$y$}
\usefont{T1}{ptm}{m}{n}
\rput(8.382343,-0.8903125){$z$}
\usefont{T1}{ptm}{m}{n}
\rput(3.2823439,-1.2503124){$\omega$}
\usefont{T1}{ptm}{m}{n}
\rput(1.8123438,0.5096875){$\omega_{V}=\sigma$}
\usefont{T1}{ptm}{m}{n}
\rput(3.6923437,1.2896875){$p$}
\usefont{T1}{ptm}{m}{n}
\rput(5.432344,0.5696875){$\omega_{H}=\tau$}
\psline[linewidth=0.04cm,arrowsize=0.05291667cm 2.0,arrowlength=1.4,arrowinset=0.4]{->}(2.0009375,-0.9203125)(4.7009373,-0.9203125)
\usefont{T1}{ptm}{m}{n}
\rput(6.632344,-0.5703125){$\gamma$}
\end{pspicture}
}
\end{center}
Then, $\omega_{V}=\sigma$, $\omega_{H}=\tau$ and
$$(g\circ \mu_{D})(\gamma\otimes \omega)=(\gamma\circ \omega)_{V}=\sigma \neq 0=\mu_{B}(id_{s(\gamma)}\otimes \sigma)=\mu_{B}\circ  (g(\gamma)\otimes g(\omega))=$$
$$(\mu_{B}\circ (g\otimes g))(\gamma\otimes \omega).$$

\end{exemplo}

\begin{prop}
Let $(D,B,f,g)$  be a weak bialgebra with a weak projection and define the morphism
$$t_{B}^{D}=\phi_{D}\circ (D\otimes (\lambda_{B}\circ g))\circ \delta_{D}:D\rightarrow D.$$
This morphism is idempotent, and moreover if $p_{B}^{D}:D\rightarrow D^{B}$, $i_{B}^{D}:D^{B}\rightarrow D$ is the splitting of $t_{B}^{D}$, the following diagram

\begin{equation}
\label{equ}
\setlength{\unitlength}{1mm}
\begin{picture}(101.00,10.00)
\put(20.00,4.00){\vector(1,0){25.00}}
\put(20.00,0.00){\vector(1,0){25.00}}
\put(55.00,2.00){\vector(1,0){21.00}}
\put(32.00,7.00){\makebox(0,0)[cc]{$\phi_{D}$ }}
\put(33.00,-4.00){\makebox(0,0)[cc]{$\beta_{D} $ }} \put(65.00,5.00){\makebox(0,0)[cc]{$p_{B}^{D} $
}} \put(13.00,2.00){\makebox(0,0)[cc]{$ D\otimes B$ }}
\put(50.00,2.00){\makebox(0,0)[cc]{$ D$ }}
\put(83.00,2.00){\makebox(0,0)[cc]{$D^{B} $ }}
\end{picture}
\end{equation}

\vspace{0.5cm}

is a coequalizer diagram where $\beta_{D}=(D\otimes (\varepsilon_{D}\circ
\phi_{D}))\circ (\delta_{D}\otimes B)$.

\end{prop}

\begin{dem}
To obtain that $t_{B}^{D}=\phi_{D}\circ (D\otimes (\lambda_{B}\circ g))\circ \delta_{D}$ is idempotent compute:
\begin{itemize}

\item[ ]$ \hspace{0.38cm}t_{B}^{D}\circ t_{B}^{D} $

\item[ ]$=\phi_{D}\circ (D\otimes (\lambda_{B}\circ g))\circ (\phi_{D}\otimes \phi_{D})\circ \delta_{D\otimes B}\circ (D\otimes (\lambda_{B}\circ g))\circ \delta_{D}   $

\item[ ]$= \mu_{D}\circ ((\mu_{D}\circ (D\otimes f))\otimes (f\circ \lambda_{B}\circ \mu_{B}\circ (g\otimes B)))\circ (D\otimes c_{D,B}\otimes B)\circ  $
\item[ ]$ \hspace{0.38cm} (\delta_{D}\otimes ((\lambda_{B}\otimes \lambda_{B})\circ c_{B,B}\circ \delta_{B}))\circ (D\otimes g)\circ \delta_{D} $

\item[ ]$=\mu_{D}\circ  $
\item[ ]$ \hspace{0.38cm}(D\otimes (\mu_{D}\circ ((f\circ \lambda_{B})\otimes (f\circ \lambda_{B}\circ \mu_{B}))\circ (c_{B,B}\otimes \lambda_{B})\circ (B\otimes c_{B,B})\circ (((g\otimes g)\circ \delta_{D})\otimes g)))\circ $
\item[ ]$ \hspace{0.38cm} (\delta_{D}\otimes D)\circ \delta_{D}$

\item[ ]$= \mu_{D}\circ (D\otimes (\mu_{D}\circ (f\otimes (f\circ \lambda_{B}\circ \Pi_{B}^{L}))\circ c_{B,B}\circ (B\otimes \lambda_{B})\circ \delta_{B}\circ g))\circ \delta_{D}  $

\item[ ]$=\mu_{D}\circ (D\otimes (f\circ \mu_{B}\circ (B\otimes \Pi_{B}^{R})\circ c_{B,B}\circ (\lambda_{B}\otimes \lambda_{B})\circ \delta_{B}\circ g))\circ \delta_{D}  $

\item[ ]$=\mu_{D}\circ (D\otimes (f\circ (id_{B}\wedge \Pi_{B}^{R})\circ \lambda_{B}\circ g))\circ \delta_{D} $

\item[ ]$= t_{B}^{D}, $

\end{itemize}
where the first equality follows by (\ref{comp}), the second one by (\ref{ant-2}) and  (\ref{bmodule}), the third one by the condition of coalgebra morphism for $g$, the associativity of $\mu_{D}$ and the coassociativity of $\delta_{D}$, the fourth one by the naturality of the braiding and the condition of coalgebra morphism for $g$, the fifth one by (\ref{pi2}) and the condition of algebra morphism for $f$. Finally, the sixth one follows by (\ref{ant-2}) and the seventh one by $id_{B}\wedge \Pi_{B}^{R}=id_{B}.$

Now we want to prove that, if $p_{B}^{D}:D\rightarrow D^{B}$, $i_{B}^{D}:D^{B}\rightarrow D$ is the splitting of $t_{B}^{D}$, then  $p_{B}^{D}$ is the coequalizer of $\phi_{D}$ and $\beta_{D}$. By (\ref{iv)}) we know that
\begin{equation}
\label{betaD}
\beta_{D}=\phi_{D}\circ (D\otimes
\overline{\Pi}_{B}^{L}).
\end{equation}
First observe that,

\begin{itemize}

\item[ ]$ \hspace{0.38cm}t_{B}^{D}\circ \phi_{D} $

\item[ ]$= \mu_{D}\circ ((\mu_{D}\circ (D\otimes f))\otimes (f\circ \lambda_{B}\circ \mu_{B}\circ (g\otimes B)))\circ \delta_{D\otimes B}$

\item[ ]$=\mu_{D}\circ (D\otimes  (f\circ \mu_{B}\circ ((\mu_{B}\circ (B\otimes \lambda_{B}))\otimes \lambda_{B})\circ (B\otimes c_{B,B})\circ (c_{B,B}\otimes B)))\circ $
\item[ ]$ \hspace{0.38cm} (((D\otimes g)\circ \delta_{D})\otimes \delta_{B})$

\item[ ]$=\phi_{D}\circ (D\otimes (\mu_{B}\circ (B\otimes \lambda_{B})\circ c_{B,B}\circ (g\otimes B)))\circ (\delta_{D}\otimes \Pi_{B}^{L})  $

\item[ ]$= \phi_{D}\circ (D\otimes (\lambda_{B}\circ \mu_{B}\circ (g\otimes B)))\circ (\delta_{D}\otimes \overline{\Pi}_{B}^{L}), $

\end{itemize}
where the first equality follows by (\ref{bmodule}) and (\ref{comp}), the second one by associativity of $\mu_{B}$, the condition of algebra morphism for $f$ and
(\ref{ant-1}), the third one by the naturality of the braiding and  the fourth one by (\ref{ant-1}) and $\Pi_{B}^{L}=\lambda_{B}\circ \overline{\Pi}_{B}^{L}$ (see (\ref{pi2})).

Now as $\beta_D$ can be written as in (\ref{betaD}), we obtain:
$$t_{B}^{D}\circ \beta_{D}=t_{B}^{D}\circ \phi_{D}\circ (D\otimes
\overline{\Pi}_{B}^{L})=\phi_{D}\circ (D\otimes (\lambda_{B}\circ \mu_{B}\circ (g\otimes B)))\circ (\delta_{D}\otimes (\overline{\Pi}_{B}^{L}\circ \overline{\Pi}_{B}^{L}) )=t_{B}^{D}\circ \phi_{D}$$
and therefore $p_{B}^{D}\circ \beta_{D}=p_{B}^{D}\circ \phi_{D}$.

It just remain to prove that $p_B^D$ satisfies the universal property of a coequalizer. But note that if $r:D\rightarrow H$ is a morphism such that $r\circ \beta_{D}=r\circ \phi_{D}$ we have
$$r\circ i_{B}^{D}\circ p_{B}^{D}=r\circ t_{B}^{D}=r\circ (id_{D}\wedge \Pi_{D}^{R})=r$$
since, by (\ref{bmodule}), $f\circ \overline{\Pi}_{B}^{R}\circ \lambda_{B}\circ g=\Pi_{D}^{R}$ holds. The morphism $r\circ i_{B}^{D}$ is the unique that makes the diagram commutes as if $s:D^{B}\rightarrow H$ is such that $s\circ p_{B}^{D}=r$, when we compose with $i_{B}^{D}$ in both members of the equality we obtain
$s=r\circ i_{B}^{D}$. Thus, (\ref{equ}) is a coequalizer diagram.
\end{dem}

\begin{exemplo}
\label{tbd}
In the conditions of Example \ref{example1}, the morphism $t_{B}^{D}$ of the previous result  is $$t_{B}^{D}(\tau)=\tau\circ \tau_{V}^{-1}=\tau_{H}.$$
 Then, in this case $D^{B}=RH$, $p_{B}^{D}(\tau)=\tau_{H}$ and $i_{B}^{D}(\omega)=\omega$.
\end{exemplo}

\begin{prop}
\label{weak-proj-cocleft}
Let $(D,B,f,g)$  be a weak bialgebra with a weak projection. Then $D\twoheadrightarrow D^{B}$ is a weak $B$-cocleft coextension with cleaving morphism $g$ and left weak inverse $g^{-1}=\lambda_{B}\circ g$.
\end{prop}

\begin{dem}
From Definitions 1.8 and 2.4 of \cite{nmra1} recall that $D$ is a weak $B$-cocleft coextension for the weak entwining structure $(B,D, \psi_{RR})$ defined in (\ref{psi}), if $(D,\phi_{D}, \delta_{D})\in {\mathcal M}_B^D(\psi_{RR})$ and there exists  morphisms $h:D\rightarrow B$, called a cleaving morphism, and  $h^{-1}:D\rightarrow B$ such that:
\begin{itemize}
\item[(i)] $h$ is a morphism of right $B$-modules.
\item[(ii)]  $h^{-1}\wedge h= e_{RR}$.
\item[(iii)]  $\mu_B\circ (B\otimes h^{-1})\circ \psi_{RR} = h^{-1}\circ \beta_{D}$.
\end{itemize}
We claim that $g:D\rightarrow B$ is a cleaving morphism. Indeed, $g$ is a right $B$-module morphism, as it is a weak projection. Now define $g^{-1} = \lambda_B\circ g$ and compute:
\begin{itemize}
\item[ ]$ \hspace{0.38cm}g^{-1}\wedge g $

\item[ ]$= \mu_B\circ ((\lambda_B\circ g) \otimes g)\circ \delta_D$

\item[ ]$= \mu_B\circ (\lambda_B\otimes B)\circ \delta_B\circ g$

\item[ ]$= \Pi_B^R\circ g$

\item[ ]$= e_{RR} ,$
\end{itemize}
hence $g^{-1}$ is a left weak inverse for $g$. Now to obtain equality $\mu_B\circ (B\otimes h^{-1})\circ \psi_{RR} = h^{-1}\circ \beta_D$ compute:
\begin{itemize}
\item[ ]$ \hspace{0.38cm} \mu_B\circ (B\otimes g^{-1})\circ \psi_{RR}$

\item[ ]$= \mu_B\circ (B\otimes \lambda_B)\circ (B\otimes (g\circ \mu_D\circ (D\otimes f)))\circ (c_{D,B}\otimes B)\circ (D\otimes \delta_B)$

\item[ ]$= \mu_B\circ (B\otimes \lambda_B)\circ (B\otimes \mu_B)\circ (c_{B,B}\otimes B)\circ (g\otimes \delta_B)$

\item[ ]$= \mu_B\circ c_{B,B}\circ ((\lambda_B\circ g)\otimes (\mu_B\circ (B\otimes \lambda_B)\circ \delta_B))$

\item[ ]$= \mu_B\circ c_{B,B}\circ ((\lambda_B\circ g)\otimes \Pi_B^L)$

\item[ ]$= \mu_B\circ c_{B,B}\circ (\lambda_B\otimes \lambda_B)\circ (g\otimes \overline{\Pi}_B^L)$

\item[ ]$= \lambda_B\circ \mu_B\circ (g\otimes \overline{\Pi}_B^L)$

\item[ ]$= g^{-1}\circ \mu_D\circ (D\otimes f)\circ (D\otimes \overline{\Pi}_B^L)$

\item[ ]$= g^{-1}\circ \beta_{D}.$

\end{itemize}
Here we used that $g$ is a morphism of right $B$-modules, the naturality of the braiding, (\ref{pi2}), the anticomultiplicativity of the antipode and the equality
\begin{equation}
\label{fpibarl}
f\circ \overline{\Pi}_B^L=\overline{\Pi}_D^L\circ f.
\end{equation}

\end{dem}

\begin{exemplo}
\label{tbd1}
 As a consequence of the previous result and by Example \ref{tbd} we obtain that if $G$ is a groupoid with exact factorization $H\bowtie V$, then
 $RG\twoheadrightarrow RH$ is a weak $RV$-cocleft coextension with cleaving morphism $g(\tau)=\tau_{V}$ and left weak inverse $g^{-1}(\tau_{V})=\tau_{V}^{-1}.$
\end{exemplo}

\begin{nota}
Let $(D,B,f,g)$  be a weak bialgebra with a weak projection and let $D\twoheadrightarrow D^{B}$ be the associated weak $B$-cocleft coextension. Then by Proposition 2.1 of
\cite{nmra1} we know that $(D^{B}, \varepsilon_{D^{B}}, \delta_{D^{B}})$ is a coalgebra in ${\mathcal C}$, where $\varepsilon_{D^{B}}:D^{B}\rightarrow K$ and
$\delta_{D^{B}}:D^{B}\rightarrow D^{B}\otimes D^{B}$ are the factorizations of $\varepsilon_{D}$ and $(p_{B}^{D}\otimes p_{B}^{D})\circ \delta_{D}$ respectively, through the coequalizer $p_{B}^{D}$. Then, as a consequence $p_{B}^{D}:D\rightarrow D^{B}$ is a coalgebra morphism.
\end{nota}

The following result is the main one of this section. Here we give the weak crossed product and the weak crossed coproduct structure of a weak bialgebra with a weak projection. This theorem generalizes the results obtained for bialgebras in braided categories given by Ardizzoni, Menini and Stefan in \cite {ardi}, and the ones by Schauenburg given in \cite{SCH2}.

\begin{teorema}
\label{principal}
 If $(D,B,f,g)$  is a weak bialgebra with a weak projection, then $D$ is a weak crossed biproduct of the coalgebra $D^{B}$ defined in (\ref{equ}) with the algebra $B$. In this case the morphisms $\psi_{B}^{D^{B}}$, $\sigma_{B}^{D^{B}}$, $\chi_{B}^{D^{B}}$ and
$\tau_{B}^{D^{B}}$ are given by
\begin{equation}
\label{psi-w1}
\psi_{B}^{D^{B}}=(p_{B}^{D}\otimes g)\circ \delta_{D}\circ \mu_{D}\circ (f\otimes i_{B}^{D}),
\end{equation}
\begin{equation}
\label{psi-w2}
\sigma_{B}^{D^{B}}=(p_{B}^{D}\otimes g)\circ \delta_{D}\circ \mu_{D}\circ (i_{B}^{D}\otimes i_{B}^{D}),
\end{equation}
\begin{equation}
\label{psi-w3}
\chi_{B}^{D^{B}}=(g\otimes p_{B}^{D})\circ \delta_{D}\circ \mu_{D}\circ (i_{B}^{D}\otimes f),
\end{equation}
\begin{equation}
\label{psi-w4}
\tau_{B}^{D^{B}}=(g\otimes g)\circ \delta_{D}\circ \mu_{D}\circ (i_{B}^{D}\otimes f).
\end{equation}
Moreover, the preunit and the precounit are
\begin{equation}
\label{proj-pre-copre}
\nu = i_{D^B\otimes B}\circ \eta_D,\;\;\;\;\upsilon=\varepsilon_{D}\circ p_{D^B\otimes B},
\end{equation}

and the following identities hold
\begin{equation}
\label{betacomul}
\delta_{D^{B}\otimes B}\circ \beta_{\nu}=(\beta_{\nu}\otimes \beta_{\nu})\circ \delta_{D},
\end{equation}
\begin{equation}
\label{tau-nabla}
\tau_{B}^{D^{B}}=(\varepsilon_{D^{B}}\otimes \delta_{B})\circ \nabla_{D^{B}\otimes B},
\end{equation}
\begin{equation}
\label{nu-nabla}
\upsilon=(\varepsilon_{D^{B}}\otimes \varepsilon_{B})\circ \nabla_{D^{B}\otimes B}.
\end{equation}
\end{teorema}

\begin{nota}
Note that condition (\ref{tau-nabla}) means that the coproduct is cosmash. If we interpret this result using a classical approach to the theory of weak crossed coproducts for weak bialgebras (i.e., dual to the one developed in \cite{ana1}), we find that this morphism $\tau_B^{D^B}$ corresponds to a trivial cycle in the weak case. Moreover, when $\tau_B^{D^B}$ is as the one given in (\ref{tau-nabla}), the coproduct is expressed as $\delta_{D^B\otimes B} = (D^B\otimes \chi_B^{D^B}\otimes B)\circ (\delta_{B^D}\otimes \delta_B)\circ \nabla_{D^B\otimes B}$. If we replace $\nabla_{D^B\otimes B}$ by $id_{D^B\otimes B}$ we obtain the dual version of the weak smash product defined by Caenepeel and De Groot in \cite{caengroot}. Finally, if we consider this weak corproduct in the Hopf algebra case, it becomes the classical cosmash coproduct. This fact is coherent with the result obtained by Schauenburg in \cite{SCH2} and by Ardizzoni, Menini and Stefan in \cite{ardi}.
\end{nota}

\noindent {\bf Proof of Theorem \ref{principal}:}

\noindent Suppose that $D$ is a weak bialgebra with a weak projection.
To obtain the weak crossed biproduct structure we will check that $(iii)$ of Theorem \ref{Teo-biproduct} holds. By assumption there exists an algebra $B$ and morphisms $f:B\rightarrow D$ of algebras and $g:D\rightarrow B$ such that $g\circ f = id_B$. There also exists a coalgebra $D^{B}$ (given by the coequalizer (\ref{equ})), a morphism $i_{B}^{D}:D^{B}\rightarrow  D$ and a morphism of coalgebras $p_{B}^{D}:D\rightarrow  D^{B}$ satisfying $p_{B}^{D}\circ i_{B}^{D} = id_{D^{B}}$.

Define:
$$i_{D^{B}\otimes B} = (p_{B}^{D}\otimes g)\circ \delta_D:D\rightarrow D^{B}\otimes B$$
and
$$p_{D^{B}\otimes B} = \mu_D\circ (i_{B}^{D}\otimes f):D^{B}\otimes B\rightarrow D.$$
Observe that $p_{D^{B}\otimes B} = \phi_D\circ (i_{B}^{D}\otimes B)$ and
$$
p_{D^{B}\otimes B}\circ i_{D^{B}\otimes B} =t_{B}^{D}\wedge (f\circ g)=id_{D}\wedge (f\circ \Pi_B^R\circ g)=id_{D}\wedge\Pi_D^R=id_{D}$$
as a consequence of being $f$ of bialgebras, $g$ of coalgebras and right $B$-modules. Hence $\nabla_{D^{B}\otimes B} = i_{D^{B}\otimes B}\circ p_{D^{B}\otimes B}$ is an idempotent morphism. It also is a morphism of right $B$-modules for $\phi_{D^{B}\otimes B}=D^{B}\otimes \mu_{B}$ and left $D^{B}$-comodules for $\rho_{D^{B}\otimes B}=\delta_{D^{B}}\otimes B$. Indeed, to check that $\nabla_{D^{B}\otimes B}$ is a morphism of right $B$ modules compute:

\begin{itemize}
\item[ ]$ \hspace{0.38cm}\nabla_{D^{B}\otimes B}\circ \phi_{D^{B}\otimes B} $

\item[ ]$= (p_{B}^{D}\otimes g)\circ \delta_{D}\circ \phi_{D}\circ ((\phi_{D}\circ (i_{B}^{D}\otimes B))\otimes B)$

\item[ ]$=((p_{B}^{D}\circ \phi_{D}) \otimes (g\circ \phi_{D}))\circ \delta_{D\otimes B}\circ ((\phi_{D}\circ (i_{B}^{D}\otimes B))\otimes B) $

\item[ ]$=((p_{B}^{D}\circ \beta_{D}) \otimes (g\circ \phi_{D}))\circ \delta_{D\otimes B}\circ ((\phi_{D}\circ (i_{B}^{D}\otimes B))\otimes B)  $

\item[ ]$= (p_{B}^{D}\otimes ((\varepsilon_{D}\otimes g)\circ \delta_{D}\circ \phi_{D}))\circ ((\delta_{D}\circ \phi_{D}\circ (i_{B}^{D}\otimes B))\otimes B)$

\item[ ]$= (p_{B}^{D}\otimes ( g\circ \phi_{D}))\circ ((\delta_{D}\circ \phi_{D}\circ (i_{B}^{D}\otimes B))\otimes B)$

\item[ ]$=\phi_{D^{B}\otimes B}\circ (\nabla_{D^{B}\otimes B}\otimes B),$

\end{itemize}
where the first equality follows by the structure of right $B$-module for $D^{B}$, the second one by (\ref{comp}), the third one by (\ref{equ}), the fourth one by the coassociativity of $\delta_{D}$ and (\ref{comp}), the fifth one by the counit condition and the last one by definition.

To prove that $\nabla_{D^{B}\otimes B}$ is a left $D^{B}$-comodule morphism, first we need to show the equality
\begin{equation}
\label{delta1}
(D\otimes t_{B}^{D})\circ \delta_{D}\circ i_{B}^{D}=\delta_{D}\circ i_{B}^{D}
\end{equation}
or, equivalently,
\begin{equation}
\label{delta2}
(D\otimes t_{B}^{D})\circ \delta_{D}\circ t_{B}^{D}=\delta_{D}\circ t_{B}^{D}
\end{equation}
and the equality
\begin{equation}
\label{second-equ}
(p_{B}^{D}\otimes D)\circ \delta_{D}\circ \phi_{D}\circ (i_{B}^{D}\otimes B)=(p_{B}^{D}\otimes \phi_{D})\circ ((\delta_{D}\circ i_{B}^{D})\otimes B).
\end{equation}

The proof of (\ref{delta2}) follows by

\begin{itemize}
\item[ ]$ \hspace{0.38cm}\delta_{D}\circ t_{B}^{D}  $

\item[ ]$=\mu_{D\otimes D}\circ (\delta_{D}\otimes (\delta_{D}\circ f\circ \lambda_{D}\circ g))\circ \delta_{D}  $

\item[ ]$= \mu_{D\otimes D}\circ (D\otimes ((D\otimes c_{D,D})\circ (((D\otimes (f\circ \lambda_{B}\circ g))\circ \delta_{D})\otimes (f\circ \lambda_{B}\circ g))\circ \delta_{D}))\circ \delta_{D} $

\item[ ]$=(\mu_{D}\otimes t_{B}^{D})\circ (D\otimes (c_{D,D}\circ ((D\otimes (f\circ \lambda_{B}\circ g))\circ \delta_{D})))\circ \delta_{D} $

\item[ ]$=(\mu_{D}\otimes (t_{B}^{D}\circ t_{B}^{D}))\circ (D\otimes (c_{D,D}\circ ((D\otimes (f\circ \lambda_{B}\circ g))\circ \delta_{D})))\circ \delta_{D} $

\item[ ]$=(D\otimes t_{B}^{D})\circ\delta_{D}\circ t_{B}^{D},  $

\end{itemize}

where the first equality is a consequence of (i) of Definition \ref{weakbialg}, the second one follows by the condition of coalgebra morphism for $f$ and $g$ as well as (\ref{ant-2}), the third one by the naturality of the braiding and the fourth one by idempotent character of $t_{B}^{D}$.

Secondly, the equality (\ref{second-equ}) follows by (i) of Definition \ref{weakbialg} and (\ref{equ}) because

\begin{itemize}
\item[ ]$ \hspace{0.38cm}(p_{B}^{D}\otimes D)\circ \delta_{D}\circ \phi_{D}\circ (i_{B}^{D}\otimes B)  $

\item[ ]$=((p_{B}^{D}\circ \phi_{D})\otimes \phi_{D})\circ \delta_{D\otimes B}\circ (i_{B}^{D}\otimes B) $

\item[ ]$= ((p_{B}^{D}\circ \beta_{D})\otimes \phi_{D})\circ \delta_{D\otimes B}\circ (i_{B}^{D}\otimes B) $

\item[ ]$=(p_{B}^{D}\otimes ((\varepsilon_{D}\otimes D)\circ \delta_{D}\circ \phi_{D}))\circ ((\delta_{D}\circ i_{B}^{D})\otimes B) $

\item[ ]$= (p_{B}^{D}\otimes \phi_{D})\circ ((\delta_{D}\circ i_{B}^{D})\otimes B).$
\end{itemize}

Then, by (\ref{second-equ}) and  (\ref{delta1}) we obtain that $\nabla_{D^{B}\otimes B}$ is a morphism of left $D^B$-comodules. Indeed:
\begin{itemize}
\item[ ]$ \hspace{0.38cm}(\rho_{D^{B}\otimes B}\otimes B)\circ  \nabla_{D^{B}\otimes B} $

\item[ ]$= (p_{B}^{D}\otimes ((p_{B}^{D}\otimes g)\circ \delta_{D}))\circ \delta_{D}\circ \phi_{D}\circ (i_{B}^{D}\otimes B) $

\item[ ]$=(p_{B}^{D}\otimes ((p_{B}^{D}\otimes g)\circ \delta_{D}\circ \phi_{D}))\circ((\delta_{D}\circ i_{B}^{D})\otimes B)  $

\item[ ]$= (p_{B}^{D}\otimes ((p_{B}^{D}\otimes g)\circ \delta_{D}\circ \phi_{D}))\circ(((D\otimes t_{B}^{D})\circ \delta_{D}\circ i_{B}^{D})\otimes B)$

\item[ ]$= (D^{B}\otimes ((p_{B}^{D}\otimes g)\circ \delta_{D}\circ \phi_{D}\circ (i_{B}^{D}\otimes B)))\circ(\delta_{D^{B}}\otimes B)$

\item[ ]$=(D^{B}\otimes \nabla_{D^{B}\otimes B})\circ \rho_{D^{B}\otimes B} .$
\end{itemize}

Finally, it is clear that the splitting of  $\nabla_{D^{B}\otimes B}$ is given by the  injection $i_{D^{B}\otimes B}$, the projection $p_{D^{B}\otimes B}$ and  the image is the object $D$.
Thus  $(iii)$ of Theorem \ref{Teo-biproduct} is satisfied considering $\omega=\varpi=id_{D}$, and thus $D$ is a weak crossed product biproduct. Now  by $(ii)\Rightarrow (i)$ of the symmetric version of Theorem \ref{uni-1},
\begin{equation}
\label{proc1}
\mu_{D^{B}\otimes B}=i_{D^{B}\otimes B}\circ \mu_{D}\circ
(p_{D^{B}\otimes B}\otimes p_{D^{B}\otimes B})
\end{equation}
and, by duality,
\begin{equation}
\label{coproc1}
\delta_{D^{B}\otimes B}=(i_{D^{B}\otimes B}\otimes i_{D^{B}\otimes B})\circ \delta_{D}\circ
p_{D^{B}\otimes B}.
\end{equation}

Moreover, by (\ref{preunit-precounit})
$$
\nu = i_{D^B\otimes B}\circ \eta_D,\;\;\;\;\upsilon=\varepsilon_{D}\circ p_{D^B\otimes B}.
$$

To obtain equality (\ref{psi-w1}), we have to use the morphism $\beta_{\nu}:B\rightarrow D^B\otimes B$ given by $\beta_{\nu} = (D^B\otimes\mu_B)\circ (\nu\otimes B)$, because, by the symmetric version of  Theorem 3.11 of \cite{mra-preunit}, we have:
\begin{equation}
\label{psidb}
\psi_{B}^{D^{B}}=\mu_{D^{B}\otimes B}\circ (\beta_{\nu}\otimes D^{B}\otimes \eta_{B}).
\end{equation}
Moreover, the following equality holds:
\begin{equation}
\label{exp-betanu}
\beta_{\nu}=((p_{B}^{D}\circ f\circ \overline{\Pi}_{B}^{L})\otimes B)\circ \delta_{B}.
\end{equation}

Indeed:
\begin{itemize}
\item[ ]$ \hspace{0.38cm} \beta_{\nu} $

\item[ ]$= (p_{B}^{D}\otimes (\mu_B\circ(g\otimes B)))\circ ((\delta_D \circ \eta_D)\otimes B)$

\item[ ]$ = (p_{B}^{D}\otimes g)\circ (D\otimes \mu_D)\circ ((\delta_D\circ \eta_D)\otimes f)$

\item[ ]$= ((p_{B}^{D}\circ \overline{\Pi}_D^L)\otimes g)\circ \delta_D\circ f$

\item[ ]$ =((p_{B}^{D}\circ f\circ \overline{\Pi}_B^L)\otimes B)\circ \delta_B,$

\end{itemize}

and, as a consequence, using (\ref{psidb}), we obtain:

\begin{itemize}
\item[ ]$ \hspace{0.38cm}\psi_{B}^{D^{B}} $

\item[ ]$=(p_{B}^{D}\otimes g)\circ \delta_{D}\circ \mu_{D}\circ ((\mu_{D}\circ ((t_{B}^{D}\circ f\circ \overline{\Pi}_{B}^{L})\otimes f)\circ \delta_{B})\otimes i_{B}^{D}) $

\item[ ]$= (p_{B}^{D}\otimes g)\circ \delta_{D}\circ \mu_{D}\circ (( f\circ ((\Pi_{B}^{L}\circ \overline{\Pi}_{B}^{L})\wedge id_{B}))\otimes i_{B}^{D}) $

\item[ ]$=(p_{B}^{D}\otimes g)\circ \delta_{D}\circ \mu_{D}\circ (( f\circ (\Pi_{B}^{L}\wedge id_{B}))\otimes i_{B}^{D})   $

\item[ ]$=(p_{B}^{D}\otimes g)\circ \delta_{D}\circ \mu_{D}\circ (f\otimes i_{B}^{D}).$

\end{itemize}
where the first equality follows by the condition of algebra morphism for $f$ and (\ref{exp-betanu}), the second one from $t_{B}^{D}\circ f=f\circ \Pi_{B}^{L}$ and the condition of algebra-coalgebra morphism for $f$, the third one by (\ref{pi1}) and the fourth one by $\Pi_{B}^{L}\wedge id_{B}=id_{B}.$

Also,  in light of the symmetric version of  Theorem 3.11 of \cite{mra-preunit}, we have
\begin{equation}
\label{sigmadb}
\sigma_{B}^{D^{B}}=\mu_{D^{B}\otimes B}\circ (D^{B}\otimes \eta_{B}\otimes D^{B}\otimes \eta_{B}).
\end{equation}
Using this formula  we obtain easily (\ref{psi-w2}).  The proofs for
(\ref{psi-w3}) and (\ref{psi-w4}) are dual and we leave the details to the reader.

To prove (\ref{betacomul}), first we need to the following  identity for $\beta_{\nu}$:
\begin{equation}
\label{newbetanu}
\beta_{\nu}=(p_{B}^{D}\otimes g)\circ \delta_{D}\circ f=i_{D^{B}\otimes B}\circ f.
\end{equation}

Indeed, as a consequence of (\ref{exp-betanu}), the equality
\begin{equation}
\label{tf}
t_{B}^{D}\circ f=f\circ \Pi_{B}^{L},
\end{equation}
the identity (\ref{pi1}) and the coalgebra condition of $f$, we have

\begin{itemize}
\item[ ]$ \hspace{0.38cm}(i_{B}^{D}\otimes B)\circ \beta_{\nu} $

\item[ ]$=((t_{B}^{D}\circ f\circ \overline{\Pi}_{B}^{L})\otimes B)\circ \delta_{B} $

\item[ ]$= (( f\circ \Pi_{B}^{L} \circ \overline{\Pi}_{B}^{L})\otimes B)\circ \delta_{B} $

\item[ ]$=(( f\circ \Pi_{B}^{L} )\otimes B)\circ \delta_{B}$

\item[ ]$= (( t_{B}^{D}\circ f)\otimes B)\circ \delta_{B}$

\item[ ]$= (t_{B}^{D}\otimes g)\circ \delta_{D}\circ f$

\end{itemize}
and this yields (\ref{newbetanu}).

Applying (\ref{newbetanu}) and (\ref{coproc1}) we obtain (\ref{betacomul})  that follows by:

\begin{itemize}
\item[ ]$ \hspace{0.38cm} \delta_{D^{B}\otimes B}\circ \beta_{\nu}$

\item[ ]$= (i_{D^{B}\otimes B}\otimes i_{D^{B}\otimes B})\circ \delta_{D}\circ f$

\item[ ]$= ((i_{D^{B}\otimes B}\circ f)\otimes (i_{D^{B}\otimes B}\circ f))\circ \delta_{B}$

\item[ ]$= (\beta_{\nu}\otimes \beta_{\nu})\circ \delta_{B}.$

\end{itemize}

Finally, using the definition of $p_{D^B\otimes B}$, compute:

\begin{itemize}
\item[ ]$ \hspace{0.38cm} \tau_{B}^{D^{B}}\circ i_{D^{B}\otimes B} $

\item[ ]$= (g\otimes g)\circ \delta_{D}\circ p_{D^{B}\otimes B}\circ i_{D^{B}\otimes B}$

\item[ ]$= (g\otimes g)\circ \delta_{D}$

\item[ ]$= \delta_{B}\circ g$

\item[ ]$=(\varepsilon_{D^{B}}\otimes \delta_{B})\circ i_{D^{B}\otimes B}$

\end{itemize}
and, composing with $p_{D^{B}\otimes B}$, we obtain (\ref{tau-nabla}).

Finally the proof for (\ref{nu-nabla}) is a consequence of the identity obtained for the precounit $\upsilon$.

{$\hfill\Box$\vspace{0.25cm}}

\begin{exemplo}
\label{example3}
As we showed in Example \ref{example1}, any finite groupoid $G$ with an exact factorization $H\bowtie V$ provides an example of a weak bialgebra with a weak projection
$$(D=RG, B=RV, f, g)$$
where $f:RV\rightarrow RG$ is defined by $f(\sigma)=\sigma$ and $g:RG\rightarrow RV$ is $g(\tau)=\tau_{V}$. By Example \ref{tbd} we know that $D^{B}=RH$ and $i_{B}^{D}(\omega)=\omega$, $p_{B}^{D}(\tau)=\tau_{H}$. Using these facts and  Theorem \ref{principal}, it is possible to compute the explicit form for the morphisms involved in the  weak crossed biproduct associated to the weak projection $(D=RG, B=RV, f, g)$. Then,
\begin{equation}
\label{4.10.1}
i_{RH\otimes RV}:RG\rightarrow RH\otimes RV, \;\;\; i_{RH\otimes RV}(\tau)=\tau_{H}\otimes \tau_{V},
\end{equation}
\begin{equation}
\label{4.10.2}
p_{RH\otimes RV}:RH\otimes RV\rightarrow RG, \;\;\; p_{RH\otimes RV}(\omega\otimes \sigma)=\left\{ \begin{array}{ll}
0\;\;\;\;\;& {\rm if}\;  s(\omega)\neq t(\sigma)
\;\\
\omega\circ \sigma &{\rm if}\;  s(\omega)=t(\sigma) \end{array} \right.,
\end{equation}
\begin{equation}
\label{4.10.3}
\nabla_{RH\otimes RV}:RH\otimes RV\rightarrow RH\otimes RV, \;\;\;
\end{equation}
$$\nabla_{RH\otimes RV}(\omega\otimes \sigma)=\left\{ \begin{array}{ll}
0\;\;\;\;\; &{\rm if}\;  s(\omega)\neq t(\sigma)
\;\\
\omega\otimes \sigma &{\rm if}\;  s(\omega)=t(\sigma) \end{array} \right.,
$$
\begin{equation}
\label{4.10.4}
\nu:R\rightarrow RH\otimes RV, \;\;\; \\ \nu(1_{R})=\sum_{x\in G_{0}}id_{x}\otimes id_{x},
\end{equation}
\begin{equation}
\label{4.10.5}
\beta_{\nu}:RV\rightarrow RH\otimes RV, \;\;\; \\ \beta_{\nu}(\sigma)=id_{t(\sigma)}\otimes \sigma,
\end{equation}
\begin{equation}
\label{4.10.6}
\upsilon: RH\otimes RV\rightarrow R, \;\;\; \\ \upsilon(\omega\otimes \sigma)=\left\{ \begin{array}{ll}
0\;\;\;\;\;& {\rm if}\;  s(\omega)\neq t(\sigma)
\;\\
1_{R} & {\rm if}\;  s(\omega)=t(\sigma) \end{array} \right.,
\end{equation}
\begin{equation}
\label{4.10.7}
\gamma_{\upsilon}:RH\otimes RV\rightarrow RH, \;\;\; \\ \gamma_{\upsilon}(\omega\otimes \sigma)=\left\{ \begin{array}{ll}
0\;\;\;\;\;& {\rm if}\;  s(\omega)\neq t(\sigma)
\;\\
\omega & {\rm if}\;  s(\omega)=t(\sigma) \end{array} \right.,
\end{equation}
\begin{equation}
\label{4.10.8}
\psi_{RV}^{RH}:RV\otimes RH\rightarrow RH\otimes RV, \;\;\; \\
\end{equation}
$$
 \psi_{RV}^{RH}(\sigma\otimes \omega)=\left\{ \begin{array}{ll}
0\;\;\;\;\;& {\rm if}\;  s(\sigma)\neq t(\omega)
\;\\
(\sigma\circ\omega)_{H}\otimes (\sigma\circ\omega)_{V} & {\rm if}\;  s(\sigma)= t(\omega) \end{array} \right.,
$$
\begin{equation}
\label{4.10.9}
\sigma_{RV}^{RH}:RH\otimes RH\rightarrow RH\otimes RV, \;\;\; \\
\end{equation}
$$
 \sigma_{RV}^{RH}(\omega^{\prime}\otimes \omega)=\left\{ \begin{array}{ll}
0\;\;\;\;\;& {\rm if}\;  s(\omega^{\prime})\neq t(\omega)
\;\\
\omega^{\prime}\circ \omega\otimes id_{s(\omega)} & {\rm if}\;  s(\omega^{\prime})= t(\omega) \end{array} \right.,
$$
\begin{equation}
\label{4.10.10}
\chi_{RV}^{RH}:RH\otimes RV\rightarrow RV\otimes RH, \;\;\; \\
\end{equation}
$$
 \chi_{RV}^{RH}(\omega\otimes \sigma)=\left\{ \begin{array}{ll}
0\;\;\;\;\;& {\rm if}\;  s(\omega)\neq t(\sigma)
\;\\
\sigma\otimes \omega & {\rm if}\;  s(\omega)= t(\sigma) \end{array} \right.,
$$
\begin{equation}
\label{4.10.11}
\tau_{RV}^{RH}:RH\otimes RV\rightarrow RV\otimes RV, \;\;\; \\
\end{equation}
$$
 \tau_{RV}^{RH}(\omega\otimes \sigma)=\left\{ \begin{array}{ll}
0\;\;\;\;\;& {\rm if}\;  s(\omega)\neq t(\sigma)
\;\\
(\sigma\circ \omega)_{V} \otimes (\sigma\circ \omega)_{V}& {\rm if}\;  s(\omega)= t(\sigma) \end{array} \right.,
$$
\begin{equation}
\label{4.10.12}
\mu_{RH\otimes RV}:RH\otimes RV\otimes RH\otimes RV\rightarrow RH\otimes RV, \;\;\; \\
\end{equation}
$$
 \mu_{RH\otimes RV}(\omega^{\prime}\otimes \sigma^{\prime}\otimes \omega\otimes \sigma)=
 $$
 $$\left\{ \begin{array}{ll}
0\;\;\;\;\;& {\rm if}\; s(\omega)\neq t(\sigma) \; {\rm or}\;  s(\omega)\neq t(\sigma) \; {\rm or}\; s(\sigma^{\prime})\neq t(\omega)
\;\\
\omega^{\prime}\circ (\sigma^{\prime}\circ \omega\circ \sigma)_{H} \otimes (\omega^{\prime}\circ \sigma^{\prime}\circ \omega)_{V}\circ \sigma & {\rm if}\;  s(\omega)= t(\sigma),\;  s(\omega)= t(\sigma), \; s(\sigma^{\prime})= t(\omega) \end{array} \right.,
$$
\begin{equation}
\label{4.10.13}
\delta_{RH\otimes RV}:RH\otimes RV\rightarrow RH\otimes RV\otimes RH\otimes RV, \;\;\; \\
\end{equation}
$$
 \delta_{RH\otimes RV}( \omega\otimes \sigma)=
\left\{ \begin{array}{ll}
0\;\;\;\;\;& {\rm if}\; s(\omega)\neq t(\sigma)
\;\\
\omega\otimes \sigma\otimes \omega\otimes \sigma & {\rm if}\;  s(\omega)= t(\sigma) \end{array} \right..
$$
\end{exemplo}

The product  obtained in Theorem \ref{principal} is an usual weak crossed biproduct together with some extra coditions, given by (\ref{betacomul}), (\ref{tau-nabla}) and (\ref{nu-nabla}). In the following result we show that these conditions are also sufficient for a weak crossed biproduct to be induced by a weak projection of weak bialgebras.

\begin{teorema}
\label{principal2}
Let $D$ be a weak bialgebra in ${\mathcal C}$ such that is a weak crossed biproduct of a coalgebra $C$ and a weak Hopf algebra $B$ with preunit $\nu$ and precounit $\upsilon$. Then, if this weak crossed biproduct satisfies
\begin{equation}
\label{betacomu2}
\delta_{C\otimes B}\circ \beta_{\nu}=(\beta_{\nu}\otimes \beta_{\nu})\circ \delta_{D},
\end{equation}
\begin{equation}
\label{tau-nabla2}
\tau_{B}^{C}=(\varepsilon_{C}\otimes \delta_{B})\circ \nabla_{C\otimes B},
\end{equation}
\begin{equation}
\label{nu-nabla2}
\upsilon=(\varepsilon_{C}\otimes \varepsilon_{B})\circ \nabla_{C\otimes B},
\end{equation}
there exists a weak projection $(D,B,f,g)$ for $D$.
\end{teorema}

\begin{dem} Let $D$ be a weak bialgebra such that  is a weak crossed biproduct of a coalgebra $C$ and a weak Hopf algebra $B$ with preunit $\nu:K\rightarrow C\otimes B$, precounit $\upsilon:C\otimes B\rightarrow K$ and associated idempotent $\nabla_{C\otimes B}$. Let $C\times B$ the image of $\nabla_{C\otimes B}$ and $\alpha:C\times B\rightarrow D$ the isomorphism of algebras and coalgebras. Then we define the morphisms $f:B\rightarrow D$ as $f=i_{B}$ and $g:D\rightarrow B$ by $g=p_{B}$ where $i_{B}$ and $p_{B}$ are the morphisms defined in (\ref{ipip}), that is
$$f=\alpha\circ p_{C\otimes B}\circ \beta_{\nu},\;\;\;\; g=(\varepsilon_{C}\otimes B)\circ i_{C\otimes B}\circ \alpha^{-1}.
$$
By (iii-1) of Theorem \ref{Teo-biproduct} we know that $f$ is an algebra morphism and $g\circ f=id_{B}$. By (\ref{nu-nabla2}) we have
$$\varepsilon_{B}\circ g=(\varepsilon_{C}\otimes \varepsilon_{B})\circ i_{C\otimes B}\circ \alpha^{-1}=\upsilon\circ i_{C\otimes B}\circ \alpha^{-1}=\varepsilon_{C\times B}\circ \alpha^{-1}=\varepsilon_{D},$$
and as a consequence $\varepsilon_{B}=\varepsilon_{B}\circ g\circ f=\varepsilon_{D}\circ f$. Moreover, by (\ref{betacomu2}), $f$ is  comultiplicative. Indeed:
\begin{itemize}
\item[ ]$ \hspace{0.38cm} \delta_{D}\circ f$

\item[ ]$= \delta_{D}\circ \alpha\circ p_{C\otimes B}\circ \beta_{\nu} $

\item[ ]$= ((\alpha\circ p_{C\otimes B}\circ \beta_{\nu})\otimes (\alpha\circ p_{C\otimes B}\circ \beta_{\nu}))\circ \delta_{B} $

\item[ ]$= (f\otimes f)\circ \delta_{B},$

\end{itemize}
and thus $f$ is a coalgebra morphism.
On the other hand, $g$ is a coalgebra morphism because:

\begin{itemize}
\item[ ]$ \hspace{0.38cm} (g\otimes g)\circ \delta_{D}$

\item[ ]$=(((\varepsilon_{C}\otimes B)\circ i_{C\otimes B})\otimes ((\varepsilon_{C}\otimes B)\circ i_{C\otimes B}))\circ \delta_{C\times B}\circ \alpha^{-1}  $

\item[ ]$= (((\varepsilon_{C}\otimes B)\circ \nabla_{C\otimes B})\otimes ((\varepsilon_{C}\otimes B)\circ \nabla_{C\otimes B}))\circ (C\otimes \chi_{B}^{C}\otimes B)\circ (\delta_{C}\otimes \tau_{B}^{C})\circ $

\item[ ]$ \hspace{0.38cm}(\delta_{C}\otimes B)\circ i_{C\otimes B} \circ \alpha^{-1}$

\item[ ]$= (((B\otimes\varepsilon_{C})\circ \chi_{B}^{C})\otimes ((B\otimes\varepsilon_{C})\circ \chi_{B}^{C}))\circ (C\otimes \chi_{B}^{C}\otimes B)\circ (\delta_{C}\otimes \tau_{B}^{C})\circ $

\item[ ]$ \hspace{0.38cm}(\delta_{C}\otimes B)\circ i_{C\otimes B} \circ \alpha^{-1}$

\item[ ]$=(B\otimes B\otimes \varepsilon_{C})\circ (B\otimes \chi_{B}^{C})\circ (\chi_{B}^{C}\otimes B)\circ (C\otimes \tau_{B}^{C})\circ (\delta_{C}\otimes B) \circ i_{C\otimes B} \circ \alpha^{-1}  $

\item[ ]$=\tau_{B}^{C}\circ \nabla_{C\otimes B}\circ  i_{C\otimes B} \circ \alpha^{-1}   $

\item[ ]$=\tau_{B}^{C}\circ   i_{C\otimes B} \circ \alpha^{-1}  $

\item[ ]$= (\varepsilon_{C}\otimes \delta_{B})\circ  i_{C\otimes B} \circ \alpha^{-1} $

\item[ ]$= \delta_{B}\circ g,$

\end{itemize}
where the first equality follows by the coalgebra morphism condition for $\alpha$, the second one by the properties of $\delta_{C\times B}$ (see (\ref{co-prod-wcp})), the third one by $(\varepsilon_{C}\otimes B)\circ \nabla_{C\otimes B}=(B\otimes\varepsilon_{C})\circ \chi_{B}^{C}$, the fourth one by (\ref{co-wmeas-wcp}) and the counit properties of $\varepsilon_{C}$, the fifth one by (\ref{co-twis-wcp}), the sixth one by the properties of the idempotent $\nabla_{C\otimes B}$, the seventh one by (\ref{tau-nabla2}) and, finally, the eighth one by definition.

To finish the proof, we only need to show  that $g$ is a morphism of $B$-modules.  To obtain this last property compute:

\begin{itemize}
\item[ ]$ \hspace{0.38cm} g\circ \mu_{D}\circ (D\otimes f)$

\item[ ]$=  (\varepsilon_{C}\otimes B)\circ i_{C\otimes B}\circ \mu_{C\times B}\circ (\alpha^{-1}\otimes (p_{C\otimes B}\circ \beta_{\nu}))$

\item[ ]$=(\varepsilon_{C}\otimes B)\circ \mu_{C\otimes B}\circ ((i_{C\otimes B}\circ  \alpha^{-1})\otimes \beta_{\nu}) $

\item[ ]$=(\varepsilon_{C}\otimes \mu_{B})\circ ((\mu_{C\otimes B}\circ ((i_{C\otimes B}\circ  \alpha^{-1})\otimes \nu) )\otimes B) $

\item[ ]$=(\varepsilon_{C}\otimes \mu_{B})\circ ((\nabla_{C\otimes B}\circ i_{C\otimes B}\circ  \alpha^{-1})\otimes B)  $

\item[ ]$= \mu_{B}\circ (g\otimes B),$

\end{itemize}
where the first equality follows by the condition of algebra morphism for $\alpha^{-1}$, the second one by the definition of $\mu_{C\times B}$ and $\nabla_{C\otimes B}\circ \beta_{\nu}=\beta_{\nu}$, the third one by the associativity of the product in $B$, the fourth one by by the symmetric equality of (\ref{idempot-preunit}) and the fifth one by the properties of $\nabla_{C\otimes B}$.

\end{dem}

\section*{Acknowledgements}
The authors were supported by Xunta de Galicia
(Project: PGIDT07PXB322079PR), Ministerio de Educaci\'on
(Projects: MTM2007-62427, MTM\-2006\--14908-CO2-01),  and  FEDER.

\end{document}